\pdfoutput=1
\RequirePackage{ifpdf}
\ifpdf
\documentclass[pdftex]{sigma}
\else
\documentclass{sigma}
\fi

\numberwithin{equation}{section}
\numberwithin{theorem}{section}
\numberwithin{proposition}{section}
\numberwithin{lemma}{section}
\numberwithin{corollary}{section}
\numberwithin{definition}{section}
\numberwithin{example}{section}
\numberwithin{remark}{section}
\numberwithin{note}{section}

\DeclareMathOperator{\id}{id}

\newcommand{\ME}[5]{E^{#1}\left(\renewcommand{\arraystretch}{0.8}
\begin{array}{cccccccc}#2
\end{array}
\renewcommand{\arraystretch}{1.0} \Big|\,{#3};{#4};{#5} \right)}

\newcommand{\SP}[3]{{{}_4 F^{n}_3} \left(\arraycolsep=1pt \renewcommand{\arraystretch}{1.2}
\begin{array}{cccccccc}#1
\end{array}
\left| \;\renewcommand{\arraystretch}{1.2}
\begin{array}{cccccccc}#2
\end{array}
\right.
\left| \;\renewcommand{\arraystretch}{1.2}
\begin{array}{cccccccc}#3
\end{array}
\;\right| \;\renewcommand{\arraystretch}{1.5}{1} \right)}

\newcommand{\SPa}[2]{{{}_4 F^{n}_3} \left(\arraycolsep=1pt \renewcommand{\arraystretch}{1.2}
\begin{array}{cccccccc}#1
\end{array}
\left| \;\renewcommand{\arraystretch}{1.2}
\begin{array}{cccccccc}#2
\end{array}
\right.
\right.\;\renewcommand{\arraystretch}{1.2}}

\begin{document}

\allowdisplaybreaks

\renewcommand{\thefootnote}{$\star$}

\renewcommand{\PaperNumber}{026}

\FirstPageHeading

\ShortArticleName{Symmetry Groups of $A_n$ Hypergeometric Series}

\ArticleName{Symmetry Groups of $\boldsymbol{A_n}$ Hypergeometric Series\footnote{This paper is a~contribution to the
Special Issue in honor of Anatol Kirillov and Tetsuji Miwa.
The full collection is available at \href{http://www.emis.de/journals/SIGMA/InfiniteAnalysis2013.html}
{http://www.emis.de/journals/SIGMA/InfiniteAnalysis2013.html}}}

\Author{Yasushi KAJIHARA}

\AuthorNameForHeading{Y.~Kajihara}

\Address{Department of Mathematics, Kobe University, Rokko-dai, Kobe 657-8501, Japan}
\Email{\href{mailto:kajihara@math.kobe-u.ac.jp}{kajihara@math.kobe-u.ac.jp}}

\ArticleDates{Received September 30, 2013, in f\/inal form March 04, 2014; Published online March 18, 2014}

\Abstract{Structures of symmetries of transformations for Holman--Biedenharn--Louck $A_n$ hypergeometric series: $A_n$
terminating balanced ${}_4 F_3$ series and $A_n$ elliptic ${}_{10} E_9$ series are discussed.
Namely the description of the invariance groups and the classif\/ication all of possible transformations for each types of
$A_n$ hypergeometric series are given.
Among them, a~``periodic'' af\/f\/ine Coxeter group which seems to be new in the literature arises as an invariance group
for a~class of $A_n$ ${}_4 F_3$ series.}


\Classification{33C67; 20F55; 33C20; 33D67}

\rightline{\it Dedicated to Professors Anatol N.~Kirillov and Tetsuji Miwa for their 65th birthday}

\renewcommand{\thefootnote}{\arabic{footnote}}
\setcounter{footnote}{0}

\section{Introduction}

In this paper, we discuss structures of symmetries of transformations for two classes of $A_n$ hypergeometric series:
$A_n$ terminating balanced ${}_4 F_3$ series and $A_n$ elliptic ${}_{10} E_9$ series.
Namely we give descriptions of the invariance groups and classif\/ication all of possible transformations for each type
of~$A_n$ hypergeometric series group-theoretically.
Among them, a~``periodic'' af\/f\/ine Coxeter group which seems to be new in the literature arises as an invariance group
for a~class of~$A_n$ ${}_4 F_3$ series.

The hypergeometric series ${}_{r+1} F_r$ is def\/ined by
\begin{gather*}
{}_{r+1} F_r \left[
\begin{matrix}
a_0, & a_1, & a_2, & \dots, & a_r
\\
& b_1, & b_2, & \dots, & b_r
\end{matrix}
; z \right]:=
\sum\limits_{k \in {\mathbb N}}
\frac {[a_0, a_1, \dots, a_r]_k}{k ! [b_1, \dots, b_r]_k} z^k,
\end{gather*}
where $[c]_k = c (c+1) \cdots (c+ k -1)$ is Pochhammer symbol and $[d_1, \dots, d_r]_k = [d_1]_k \cdots [d_r]_k$.

Investigations of the symmetry of the hypergeometric series goes back to 19th century in the case of ${}_3 F_2$ series.
Thomae~\cite{Thomae} has considered the following ${}_3 F_2$ transformation formula
\begin{gather*}
%\label{3F2trans}
{}_3 F_2 \left[
\begin{matrix}
a, b, c
\\
d, e
\end{matrix}
; 1 \right] = \frac{\Gamma(e) \Gamma(d+e-a- b-c)}{\Gamma(e- a) \Gamma(d+e-b- c)}{}_3 F_2 \left[
\begin{matrix}
a, d-b, d- c
\\
d, d+e-b-c
\end{matrix}
; 1 \right],
\end{gather*}
where $\Gamma(x)$ is the Euler gamma function.
Later, Hardy~\cite{Hardy} formulated this case as follows, where we give a~ref\/ined form (see also
Whipple~\cite{Whipple1}):

\medskip

\noindent {\bf Theorem} (Hardy). {\it Let $s = s(x_1, x_2, x_3, x_4, x_5) = x_1+x_2+x_3-x_4-x_5$.
The function
\begin{gather*}
%\label{ThRel}
\frac{1}{\Gamma(s) \Gamma(2 x_4) \Gamma(2 x_5)}{}_3 F_2 \left[
\begin{matrix}
2 x_1 -s, 2 x_2 -s, 2 x_3 -s
\\
2 x_4, 2 x_5
\end{matrix}
; 1 \right]
\end{gather*}
is a~symmetric function of the $5$ variables $x_1$, $x_2$, $x_3$, $x_4$, $x_5 $.
Thus the ${}_3 F_2$ series have a~symmetry of the symmetric group~$\mathfrak{S}_5$ of degree~$5$.
}

 \medskip

  Motivated by quantum mechanics and representation theory, the symmetry of hypergeometric series has been
investigated by many authors including physicists.
It is also related to highest weight representations of the unitary group~${\rm SU}(2)$.
(for a~expository in this direction, we refer to the paper by Krattenthaler and Srinivasa Rao~\cite{K-SR}).
The Clebsch--Gordan coef\/f\/icients can be expressed in terms of ${}_3 F_2$ series, in other words, Hahn polynomials and by
using Hardy's result, one f\/inds non-trivial zeros of the coef\/f\/icients.
That is, one can clarify the structure of the highest weight representations.
The Racah coef\/f\/icients can be expressed in terms of terminating balanced ${}_4 F_3$ series, in other words, Racah
polynomials.
The corresponding results for the ${}_4 F_3$ series has been given by Beyer, Louck and Stein~\cite{Louck} (see also
Section~\ref{section2.2}).
The results of the groups of symmetry for hypergeometric series have been generalized for each types of hypergeometric
series (see~\cite{LvdJ, vdBRS,vdj}).

We also mention that recently, number-theorists have investigated in this direction: Formicella, Green and
Stade~\cite{Stade} and Mishev~\cite{Mishev} discussed in the case of non-terminating (but) ba\-lan\-ced ${}_4 F_3$ series
with a~connection with Fourier coef\/f\/icients of $GL_n$ automorphic form.
In~\cite{KR}, Krattenthaler and Rivoal presented a~dif\/ferent but considerably interesting approach related to their
investigations regarding odd values for Riemann zeta functions.

Elliptic hypergeometric series has f\/irst introduced by Frenkel and Turaev~\cite{FT} in the context of elliptic
$6j$-symbol.
They obtained transformation and summation formulas for elliptic hypergeometric series by using invariants of links
which extends the works by A.N.~Kirillov and N.Yu.~Reshetikhin~\cite{KirR} (see also~\cite{Kir}).

\looseness=1
In 1970's, Holman, Biedenharn and Louck~\cite{HBL} and Holman~\cite{H1} has introduced a~class of multivariate
generalization of hypergeometric series which is nowadays called as $A_n$ hypergeometric series (or hypergeometric
series in ${\rm SU}(n+1)$) for explicit expressions of Clebsch--Gordan and Racah coef\/f\/icients of the higher dimensional unitary
group ${\rm SU}(n+1)$.
It includes $A_n$ ${}_4 F_3$ series which we discuss in Section~\ref{section2}.
Results of transformation and summation formulas for $A_n$ hypergeometric series including basic and elliptic
generalization and extension to other (classical) root systems has known by many authors (for summary, see an excellent
exposition by S.C.~Milne~\cite{MilneNagoya}).

Among them, we obtained a~number of transformation formulas for (mainly basic) hypergeometric series of type $A$ with
dif\/ferent dimensions in~\cite{Kaji1} (see also~\cite{KajiS} and~\cite{KajiR}).
In the joint work with M.~Noumi~\cite{KajiNou}, we showed the results can be extended in the case of balanced series and
proposed the notion of {\it duality transformation formula}.
In~\cite{Kaji1} and~\cite{KajiNou}, we have obtained our results by starting from the Cauchy kernels and their action of
($q$-)dif\/ference operators of Macdonald type.
The class of hypergeometric transformations of type $A$ with dif\/ferent dimensions in our previous works can be
considered to involve some of previously known $A_n$ hypergeometric transformation formulas in 20th century
(see~\cite{MilneNagoya}).
In~\cite{KajiBDT} (see also~\cite{KajiS} and~\cite{KajiNou}), we proved a~number of their results by combining some
special cases (hypergeometric transformations between $A_n$ hypergeometric series and one-dimensional ($A_1$)
hypergeometric series).
This paper can be considered to be a~continuation of~\cite{KajiBDT}.

In this paper, we discuss the symmetry of some classes of $A_n$ hypergeometric series including $n=1$ case.
Namely we investigate the invariance forms and the groups describing the symmetry of each type of hypergeometric series.
For $n \ge 2$, the symmetry of the $A_n$ hypergeometric series is more restricted than $n=1$ case if we f\/ix the symmetry
corresponding to the dimension of the summation.
So, the groups of symmetry are subgroups of that in the case of $n=1$.
Furthermore, we classify all the hypergeometric transformations which can be obtained by the combinations of possible
permutations of the parameters and the hypergeometric transformations without trivial transformations in each cases.
The classif\/ications are given by double coset decomposition of the corresponding groups.

In Section~\ref{section2}, we discuss symmetries of $A_n$ terminating balanced ${}_4 F_3$ series.
Among these, a~{\it ``periodic affine''} Weyl group that is periodic with respect to the translations arises in a~class
of $A_n$ ${}_4 F_3$ series.
It seems not to have previously appeared in the literature as a~Coxeter group (see~\cite{Bour,HBook} and the
paper by Iwahori and Matsumoto~\cite{IM} regarding the Weyl groups with translations).
In Section~\ref{section3}, we discuss symmetries of $A_n$ elliptic hypergeometric series.
What is remarkable in this case is a~subgroup structure.

It would be interesting if the discussions and results does work for future works not only for multivariate
hypergeometric transformations themselves, but also for deeper investigations to the structure of irreducible
decompositions of the tensor products of certain representations of higher dimensional unitary group ${\rm SU}(n+1)$ and
elliptic quantum groups of ${\rm SU}(n+1)$, the original problem to introducing $A_n$ and elliptic hypergeometric series.

On the other hand, Kajiwara et~al.~\cite{KMNOY} found that elliptic hypergeometric series ${}_{10} E_9$ arises as a~class of solutions of the elliptic
Painlev\'{e} equation associated to the af\/f\/ine Weyl group~$W(E_7^{(1)})$ which is the one of the family of the
Painlev\'{e} equations introduced by Sakai~\cite{Sakai} from the geometry of rational surfaces.
We also mention the work of Rains~\cite{Rains0} on relations between elliptic hypergeometric integrals and tau functions
of elliptic Painlev\'{e} equations (see also~\cite{Rains1} and~\cite{vdBRS}).
It would be an interesting problem to give geometric interpretation of the symmetries of classes of $A_n$ hypergeometric
series in terms of certain rational surfaces.

\section[Symmetry groups of $A_n$ ${}_4 F_3$ series]{Symmetry groups of $\boldsymbol{A_n}$ $\boldsymbol{{}_4 F_3}$ series}
\label{section2}

\subsection[Preliminaries on $A_n$ hypergeometric series]{Preliminaries on $\boldsymbol{A_n}$ hypergeometric series}

Here, we note the conventions for naming series as~$A_n$ (ordinary) hypergeometric series (or hypergeometric series in
${\rm SU}(n+1)$).
Let $\gamma= (\gamma_1, \dots, \gamma_n) \in \mathbb{N}^n$ be a~multi-index.
We denote
\begin{gather*}
\Delta(x):= \prod\limits_{1 \le i < j \le n} (x_i-x_j)
\qquad
\mbox{and}
\qquad
\Delta(x+\gamma):= \prod\limits_{1 \le i < j \le n} (x_i+{\gamma_i}-x_j-{\gamma_j}),
\end{gather*}
as the Vandermonde determinant for the sets of variables $x=(x_1, \dots, x_n)$ and $x+\gamma = (x_1+{\gamma_1},
\dots, x_n +{\gamma_n})$ respectively.
In this paper we refer multiple series of the form
\begin{gather}
\label{DefAnHS}
\sum\limits_{\gamma \in {\mathbb N}^n} \frac{\Delta (x +{\gamma})}{\Delta (x)} H(\gamma)
\end{gather}
which reduce to hypergeometric series ${}_{r+1} F_r$ for a~nonnegative integer $r$ when $n=1$ and symmetric with respect
to the subscript $1 \le i \le n$ as~$A_n$ hypergeometric series.
We call such a~series balanced if it reduces to a~balanced series when $n=1$.
Terminating, balanced and so on are def\/ined similarly.
The subscript $n$ in the label $A_n$ attached to the series is the dimension of the multiple series~\eqref{DefAnHS}.

Before beginning our discussion, we summarize $q \to 1$ results of $A_n$ Sears transformation from~\cite{KajiBDT} which
we discuss in this paper.
For the procedure of $q \to 1$ limit, one can f\/ind in the book by Gasper--Rahman~\cite{GR1} (see also~\cite{Kaji1}).

We introduce the notation for $A_n$ ${}_4 F_3$ series as follows
\begin{gather*}
%\label{defPhiSer}
\SP{\{b_i\}_n
\\
\{x_i\}_n}{a_1, a_2
\\
e_1, e_2}{c\\d}
\\
\qquad
:=
\sum\limits_{\gamma \in {\mathbb N}
^n} \frac{\Delta (x+{\gamma})}{\Delta (x)} \prod\limits_{1 \le i, j \le n} \frac{[ b_j+x_i-x_j]_{\gamma_i}}
{[ 1+x_i-x_j]_{\gamma_i}}
\prod\limits_{1 \le i \le n} \left[ \frac{[ c+x_i ]_{\gamma_i}}{[d+x_i ]_{\gamma_i}} \right]
\frac{[ a_1, a_2]_{|\gamma|}}{[e_1, e_2 ]_{|\gamma|}},
\end{gather*}
where $|\gamma| = \gamma_1+\gamma_2+\dots+\gamma_n$ is a~length of the multi-index $\gamma$.

Here we give two $A_n$ Whipple transformations which discuss in this paper.

 {\bf Rectangular version} (the ${q \to 1}$ limit of ${A_n}$
Sears transformation formula, Corollary~4.5 in~\cite{KajiBDT})
\begin{gather}
\SP{\{{-m_i}\}_n\\ \{x_i\}_n}{a_1, a_2\\e_1, e_2}{c\\d}
= \frac{[ d+e_1-a_2-c]_{|M|}}{[ d+e-a_1-a_2-c]_{|M|}} \prod\limits_{1 \le i \le n} \frac{[ d-a_1+x_i]_{m_i}}{[d+x_i]_{m_i}}
\nonumber
\\
\qquad
\times \SP{\{{-m_i}\}_n
\\
\{\tilde{x}_i\}_n}{a_1, e_1-a_2
\\
e_1, d+e_1-a_2-c}{e_1-c
\\
e_1+e_2-a_2-c},
\label{ias2}
\end{gather}
where $|M| = m_1+m_2+\cdots+m_n$ and $\tilde{x}_i ={-m_i+|M|}-x_i$ for $1\le i \le n$.
The balancing condition in this case is
\begin{gather*}
a_1+a_2+c+{1-|M|} = d+e_1+e_2.
\end{gather*}

Note that $\SP{\{{-m_i}\}_n
\\ \{x_i\}_n}{a_1, a_2
\\e_1, e_2}{c\\d}$ series terminates with respect to a~multi-index.
In this paper, we call such series as~{\it rectangular} and the multiple series which terminates with respect to the
length of multi-indices as~{\it triangular}.

  {\bf Triangular version} (the ${q \to 1}$ limit of ${A_n}$
Sears transformation formula, Proposition~4.5
in~\cite{KajiBDT})
\begin{gather}
\SP {\{b_i\}_n\\ \{x_i\}_n}{{-N}, a \\e_1, e_2}{c\\d}
= \frac{[d+e_1-a - c]_N}{[ d+e_1-a - B-c]_N} \prod\limits_{1 \le i \le n} \frac{[ d-b_i+x_i]_N}{[d+x_i]_N}
\nonumber
\\
\qquad{}
\times \SP {\{b_i\}_n
\\
\{\tilde{x}_i\}_n}{{-N}, e_1-a\\
e_1, d+e_1-a - c}{e_1-c
\\
e_1+e_2-a - c},
\label{las2}
\end{gather}
where $|B|= b_1+b_2+\dots+b_n$ and $\tilde{x}_i = b_i-B-x_i$ for $1 \le i \le n$.
The balancing condition in this case is
\begin{gather*}
a~+ B+c+{1-N} = d+e_1+e_2.
\end{gather*}

\begin{remark}
In the case when $n=1$ and $x_1 = 0$,~\eqref{ias2} and~\eqref{las2} reduce to the Whipple transformation formula for
terminating balanced ${}_4 F_3$ series
\begin{gather}
{}_{4} F_{3} \left[
\begin{matrix}
{-N}, a_1, a_2, a_3
\\
d_1, d_2, d_3
\end{matrix}
; 1 \right] = \frac{[ d_2-a_1, d_1+d_2-a_2-a_3]_N}{[d_2, d_1+d_2-a_1-a_2-a_3]_N}
\nonumber
\\
\hphantom{{}_{4} F_{3} \left[
\begin{matrix}
{-N}, a_1, a_2, a_3
\\
d_1, d_2, d_3
\end{matrix}
; 1 \right] =}{}
\times {}_4 F_3 \left[
\begin{matrix}
{-N}, a_1, d_1-a_3, d_1-a_2
\\
d_1, d_1+d_3-a_2-a_3, d_1+d_2-a_2-a_3
\end{matrix}
; 1 \right].
\label{d1st1}
\end{gather}
Note that identity above~\eqref{d1st1} is valid if the balancing condition
\begin{gather*}
a_1+a_2+a_3+{1-N} = d_1+d_2+d_3
\end{gather*}
holds.

\end{remark}

\subsection[Symmetries of ${}_4 F_3$ transformations ($A_1$ case)]{Symmetries
of $\boldsymbol{{}_4 F_3}$ transformations ($\boldsymbol{A_1}$ case)}
\label{section2.2}

Here,
we discuss the symmetry for terminating balanced ${}_4 F_3$ series, namely the $A_1$ case:
\begin{gather}
\label{d1Ser}
{}_{4} F_{3} \left[
\begin{matrix}
{-N}, a_1, a_2, a_3
\\
d_1, d_2, d_3
\end{matrix}
; 1 \right]
\end{gather}
with the balancing condition
\begin{gather}
\label{d1bc}
a_1+a_2+a_3+{1-N} = d_1+d_2+d_3.
\end{gather}
Though most of all the results were originally obtained in~\cite{Louck} with a~dif\/ferent formulation, we continue our
discussion in order to give a~correspondence to the results in Section~\ref{section2.5}.

The action of the parameters $a_i$, $d_i$, $i= 1, 2, 3$, for the Whipple transformation~\eqref{d1st1} is given as~follows
\begin{gather*}
%\label{actd1s}
s: \
\begin{bmatrix}
a_1^{}
\\
a_2^{}
\\
a_3
\\
d_1
\\
d_2^{}
\\
d_3^{}
\end{bmatrix}
\to
\begin{bmatrix}
a_1
\\
d_1-a_3
\\
d_1-a_2
\\
d_1
\\
d_1+d_2-a_2-a_3
\\
d_1+d_3-a_2-a_3
\end{bmatrix}
.
\end{gather*}
One can consider it as a~linear transformation acting on the vector $\vec{v}_1 = {}^t (a_1, a_2, a_3, d_1, d_2, d_3)$.
The matrix realization $S$ for transformation $s$ is given as follows
\begin{gather*}
%\label{Matrd1S1}
S:=
\begin{bmatrix}
1 & 0 & 0 & 0 & 0 & 0
\\
0 & 0 & -1 & 1 & 0 & 0
\\
0 & -1& 0 & 1 & 0 & 0
\\
0 & 0 & 0 & 1 & 0 & 0
\\
0 & -1& -1 & 1 & 1 & 0
\\
0 & -1& -1 & 1 & 0 & 1
\end{bmatrix}
.
\end{gather*}

It is easy to see that the ${}_4 F_3$ series is invariant under the action of the permutation in the two sets of
parameters $\{a_1, a_2, a_3\}$ and $\{d_1, d_2, d_3\}$.
For $i=1, 2$, let $r_i$ be the permutation of~$a_i$ and~$a_{i+1}$ and let~$t_i$ be the permutation of~$d_i$ and~$d_{i+1}$.
The matrix realizations $R_i$ (resp.~$T_i$) of~$r_i$ (resp.~$t_i$) is given by its action on the vector~$\vec{v}_1 $.
For example,
\begin{gather*}
%\label{MatrPerm}
R_1:=
\begin{bmatrix}
0 & 1 & 0 & 0 & 0 & 0
\\
1 & 0 & 0 & 0 & 0 & 0
\\
0 & 0 & 1 & 0 & 0 & 0
\\
0 & 0 & 0 & 1 & 0 & 0
\\
0 & 0 & 0 & 0 & 1 & 0
\\
0 & 0 & 0 & 0 & 0 & 1
\end{bmatrix}
\qquad
\text{and}
\qquad
T_1:=
\begin{bmatrix}
1 & 0 & 0 & 0 & 0 & 0
\\
0 & 1 & 0 & 0 & 0 & 0
\\
0 & 0 & 1 & 0 & 0 & 0
\\
0 & 0 & 0 & 0 & 1 & 0
\\
0 & 0 & 0 & 1 & 0 & 0
\\
0 & 0 & 0 & 0 & 0 & 1
\end{bmatrix}
.
\end{gather*}

We have the invariant form for the terminating balanced ${}_4 F_3$ series:

\begin{proposition}[invariance for terminating balanced ${{}_4 F_3}$ series]
\begin{gather*}
%\label{inv4F3}
{}_4 \tilde{F}_3 \left[ \vec{v}_1 \right] := [ d_1, d_2, d_3]_N {}_{4} F_{3} \left[
\begin{matrix}
{-N}, a_1, a_2, a_3
\\
d_1, d_2, d_3
\end{matrix}
; 1 \right]
\end{gather*}
is invariant under all of the actions $r_i$, $t_i$, $i=1, 2$, and $s$.
\end{proposition}

Obviously, the transformations $r_1$ and $r_2$ enjoy the braid relation $r_1 r_2 r_1 = r_2 r_1 r_2$ and $r_i^2 = \id$
for $i=1, 2$.
So do $t_1$ and $t_2$.
The relations among the element $s$ and others are summarized as follows:

\begin{lemma}
We have $(s t_1)^3 = (s r_1)^3 =\id
$ and $(s t_2)^2 = (s r_2)^2 =\id$.
\end{lemma}

We def\/ine the mapping $\pi_1$ as
\begin{gather*}
%\label{TaioPi1}
s_1 \to
\sigma_3,
\qquad
 r_i \to \sigma_{3-i}
,
\qquad
 t_i \to
\sigma_{3+i},
\qquad
i= 1, 2.
\end{gather*}
Then, by braid relations among $r_i$ and $t_i$ and the lemma above, we see that the following relation holds:
\begin{gather*}
\begin{cases}
\sigma_i \not= \id,
\qquad
\sigma_i^2 = \id, &  i = 1, 2, 3, 4, 5 ,
\\
\sigma_{i}
\sigma_{i+1}
\sigma_{i}
=
\sigma_{i+1}
\sigma_{i}
\sigma_{i+1},
&  i = 1, 2, 3, 4 ,
\\
\sigma_{i}
\sigma_{j}
=
\sigma_{j}
\sigma_{i},
&  |i-j| \ge 2 .
\end{cases}
\end{gather*}
Thus we have the following.

\begin{proposition}
The set generated by $r_i$, $t_i$ for $i=1, 2$ and $s$ forms a~Coxeter group.
Furthermore, the group is isomorphic to $\mathfrak{S}_6$.
\end{proposition}

Here we classify the possible transformations for terminating balanced ${}_4 F_3$ series of the form~\eqref{d1Ser}.
Recall that ${}_4 F_3$ series is invariant under the action $\sigma_k$ for $k= 1, 2, 4, 5$.
Thus our problem reduces to give an orbit decomposition of the double coset $H \backslash G /H$,
where $G:= \{\sigma_i \,|\, i= 1, 2, 3,4, 5\}$, $G_1:= \{\sigma_i \,|\, i= 1, 2\}$,
$G_2:= \{\sigma_i \,|\, i= 4, 5\}$ and $H = G_1 \times G_2$.
The representatives of orbits of $H \backslash G / H$ is given by
\begin{alignat*}{3}
&(i)\quad &&\omega_0 =\id,&\\
&(ii)\quad &&\omega_1=\sigma_3,&\\
& (iii)\quad &&\omega_2= \sigma_3 \sigma_4 \sigma_2 \omega_1 = \sigma_3 \sigma_4 \sigma_2 \sigma_3,&\\
& (iv)\quad && \omega_3 = \sigma_3 \sigma_4 \sigma_5 \sigma_2 \sigma_1 \omega_2 = \sigma_3 \sigma_4 \sigma_5 \sigma_2 \sigma_1
\sigma_3 \sigma_4 \sigma_2 \sigma_3.&
\end{alignat*}

Thus we are ready to present a~list of the ${}_4 F_3$ transformations according to the representatives above.
We frequently make the simplif\/ication of the product factor by using the balancing condition~\eqref{d1bc}.

The transformation associated with $(i)$ is identical.
The second one $(ii)$ is the Whipple transformation~\eqref{d1st1} itself.
The third one $(iii)$ is given by
\begin{gather}
{}_{4} F_{3} \left[
\begin{matrix}
{-N}, a_1, a_2, a_3
\\
d_1, d_2, d_3
\end{matrix}
; 1 \right]
= \frac{[d_1+d_2-a_1-a_3, d_1+d_2-a_2-a_3, a_3]_N}{[d_1, d_2, d_1+d_2-a_1-a_2-a_3]_N}
\nonumber
\\
\qquad
\times {}_4 F_3 \left[
\begin{matrix}
{-N}, d_1-a_3, d_2-a_3, d_1+d_2-a_1-a_2-a_3
\\
d_1+d_2-a_2-a_3, d_1+d_2-a_1-a_3, d_1+d_2+d_3-a_1-a_2-2 a_3
\end{matrix}
; 1 \right].
\label{d1rst1}
\end{gather}
The forth one $(iv)$ is
\begin{gather}
{}_{4} F_{3} \left[
\begin{matrix}
{-N}, a_1, a_2, a_3
\\
d_1, d_2, d_3
\end{matrix}
; 1 \right] = \frac{[a_1, a_2, a_3]_N}{[d_1, d_2, d_1+d_2-a_1-a_2-a_3 ]_N}
\nonumber
\\
\qquad\phantom{=}\times {}_4 F_3 \left[
\begin{matrix}
{-N}, d_1+d_2-a_1-a_2-a_3, d_1+d_3-a_1-a_2-a_3,
\\
d_1+d_2+d_3-a_1-a_2-2 a_3, d_1+d_2+d_3 -
a_1-2 a_2-a_3,
\end{matrix}
\right.
\nonumber
\\
\hspace*{45mm}
\left.
\begin{matrix}
d_2+d_3-a_1-a_2-a_3
\\
d_1+d_2+d_3-2 a_1-a_2-a_3
\end{matrix}
; 1 \right].
\label{d1r1}
\end{gather}

The transformation~\eqref{d1r1} has an alternative expression
\begin{gather}
{}_{4} F_{3} \left[
\begin{matrix}
{-N}, a_1, a_2, a_3
\\
d_1, d_2, d_3
\end{matrix}
; 1 \right] = \left(-1 \right)^N \frac{[a_1, a_2, a_3]_N}{[d_1, d_2, d_3]_N}
\nonumber\\
 \phantom{{}_{4} F_{3} \left[
\begin{matrix}
{-N}, a_1, a_2, a_3
\\
d_1, d_2, d_3
\end{matrix}
; 1 \right] =}\times {}_4 F_3 \left[
\begin{matrix}
{-N}, 1-N-d_1, 1-N-d_2, 1-N-d_3
\\
1-N-a_1, 1- N-a_2, 1-N-a_3
\end{matrix}
; 1 \right].
\label{ad1r1}
\end{gather}

Note also that~\eqref{ad1r1} is an inversion of the order of the summation in the ${}_4 F_3$ series.

\subsection[Symmetry of $A_n$ ${}_4 F_3$ series of rectangular type]{Symmetry
of $\boldsymbol{A_n}$ $\boldsymbol{{}_4 F_3}$ series of rectangular type}
\label{section2.3}

Here, we describe the invariance group for $A_n$ Whipple transformation of rectangular type~\eqref{ias2}, namely the
series of the form
\begin{gather}
\label{itSer}
\SP{\{{-m_i}\}_n
\\
\{x_i\}_n}{a_1, a_2
\\
e_1, e_2}{c\\d}
\end{gather}
with the balancing condition
\begin{gather}
\label{ibc}
a_1+a_2+ c+{1-|M|} = d+e_1+e_2.
\end{gather}

Suppose that $n \ge 2$ till stated otherwise.
Hereafter, we also f\/ix the symmetry of the dimension of the summation in the multiple series.

Recall that the transformation of coordinates in the right hand side of~\eqref{ias2} is given as follows
\begin{gather*}
%\label{actS}
s: \
\begin{bmatrix}
a_1^{}
\\
a_2^{}
\\
c
\\
d
\\
e_1^{}
\\
e_2^{}
\end{bmatrix}
\mapsto
\begin{bmatrix}
a_1^{1}
\\
a_2^{1}
\\
c^{1}
\\
d^{1}
\\
e_1^{1}
\\
e_2^{1}
\end{bmatrix}
=
\begin{bmatrix}
a_1
\\
e_1-a_2
\\
e_1-c
\\
e_1+e_2-a_2-c
\\
e_1
\\
e_1+d -a_2-c
\end{bmatrix}
.
\end{gather*}
It is easy to see that this transformation of coordinates is linear for $a_1$, $a_2$, $c$, $d$, $e_1$ and $e_2$.
Thus we give a~$6 \times 6$ matrix realization for transformation for $s$ acting on the vector $\vec{v} = {}^t [a_1,
a_2, c, d, e_1, e_2]$ as follows
\begin{gather*}
%\label{MatrS}
\begin{bmatrix}
a_1^{1}
\\
a_2^{1}
\\
c^{1}
\\
d^{1}
\\
e_1^{1}
\\
e_2^{1}
\end{bmatrix}
=S_1 \vec{v},
\qquad
S_1:=
\begin{bmatrix}
1 & 0 & 0 & 0 & 0 & 0
\\
0 & -1 & 0 & 0 & 1 & 0
\\
0 & 0 & -1 & 0 & 1 & 0
\\
0 & -1 & -1 & 0 & 1 & 1
\\
0 & 0& 0 & 0 & 1& 0
\\
0 & -1 & -1 & 1 & 1& 0
\end{bmatrix}
.
\end{gather*}
Note that the series~\eqref{itSer} is symmetric with respect to the two sets of parameters $\{a_1, a_2\}$ and $\{e_1,
e_2\}$.
Let $s_0$ be a~permutation of $a_1$ and $a_2$ and let $s_2$ be a~permutation of $e_1$ and $e_2$.
The matrix realization $S_0$ (resp.~$S_2$) of $s_0$ (resp.~$s_2$) is given by
\begin{gather*}
%\label{MatrR}
S_0:=
\begin{bmatrix}
0 & 1 & 0 & 0 & 0 & 0
\\
1 & 0 & 0 & 0 & 0 & 0
\\
0 & 0 & 1 & 0 & 0 & 0
\\
0 & 0 & 0 & 1 & 0 & 0
\\
0 & 0 & 0 & 0 & 1 & 0
\\
0 & 0 & 0 & 0 & 0 & 1
\end{bmatrix}
\qquad
\text{and}
\qquad
S_2:=
\begin{bmatrix}
1 & 0 & 0 & 0 & 0 & 0
\\
0 & 1 & 0 & 0 & 0 & 0
\\
0 & 0 & 1 & 0 & 0 & 0
\\
0 & 0 & 0 & 1 & 0 & 0
\\
0 & 0 & 0 & 0 & 0 & 1
\\
0 & 0 & 0 & 0 & 1 & 0
\end{bmatrix}
.
\end{gather*}
For the action on the variables $x_i$, we set to be $s_1 \cdot x_i = \tilde{x}_i = -m_i+|M|-x_i$ and otherwise to be
identical.

We introduce the normalized ${}_4 F^n_3$ series ${}_4 \widetilde{F}^n_{3} \left((\vec{v}, x)\right)$ as
\begin{gather}
\label{R4F3}
{}_4 \widetilde{F}^n_{3} \left((\vec{v}, x) \right):= [e_1, e_2]_{|M|} \prod\limits_{1 \le i \le n} [ d+x_i]_{m_i}
\SP{\{{-m_i}\}_n
\\
\{x_i\}_n}{a_1, a_2
\\
e_1, e_2}{c\\d},
\end{gather}
under the balancing condition~\eqref{ibc}.

\begin{proposition}[invariance form for multiple series of type~\eqref{itSer}] \label{proposition3}
${}_4\widetilde{F}^n_{3}\left((\vec{v},x)\right)$ is invariant under the action of~$s_0$,~$s_1$ and~$s_2$.
\end{proposition}

In the course of the proof of this proposition, we use the following lemma.

\begin{lemma}\label{lemma2}
If~\eqref{ibc} holds, then we have the following
\begin{gather*}
1)\quad [b]_{|M|} = (-1)^{|M|} [ d+e_1+e_2-a_1 -a_2-c- b]_{|M|},
\\
2)\quad[f+x_i ]_{m_i} = (-1)^{m_i} [ d+e_1+e_2-a_1 -a_2-c+ \tilde{x}_i- f]_{m_i}.
\end{gather*}
\end{lemma}

One can prove this lemma by using an elementary manipulation of shifted factorials $[z]_m = (-1)^m [1-z -m]_m$ and the
balancing condition~\eqref{ibc}.

\begin{proof}[Proof of Proposition~\ref{proposition3}]
Since the ${}_4 \widetilde{F}^n_{3} \left((\vec{v}, x)\right)$ is
symmetric with respect to the sets of parameters $\{a_1, a_2\}$ and $\{e_1, e_2\}$, it is obvious in the case of
$s_0$ and $s_2$.
For the case of $s_1$, by using the transformation formula~\eqref{ias2} and Lemma~\ref{lemma2}
\begin{gather*}
{}_4 \widetilde{F}^n_{3} \left(s_1 (\vec{v}, x) \right) = {}_4 \widetilde{F}^n_{3} \left((S_1 \vec{v}, \tilde{x}
)\right)
= [e_1, d+e_1-a_2-c]_{|M|} \prod\limits_{1 \le i \le n} [ d-a_1+x_i]_{m_i}
\\
\qquad
\phantom{=}{}
\times \SP{\{{-m_i}\}_n
\\
\{\tilde{x}_i\}_n}{a_1, e_1-a_2
\\
e_1, d+e_1-a_2-c}{c\\d}
\\
\qquad{}
= [e_1, d+e_1-a_2-c]_{|M|} \prod\limits_{1 \le i \le n} [ d-a_1+x_i]_{m_i}
\\
\qquad
\phantom{=}{}
\times \frac{[ d+e_1-a_1-a_2-c]_{|M|}}{[ d+e-a_2-c]_{|M|}} \prod\limits_{1 \le i \le n} \frac{[ d+x_i]_{m_i}}{[d-a_1+x_i]_{m_i}}
\, \SP{\{{-m_i}\}_n\\ \{x_i\}_n}{a_1, a_2\\e_1, e_2}{c\\d}
\nonumber
\\
\qquad
= {}_4 \widetilde{F}^n_{3} \left((\vec{v}, x) \right).
\end{gather*}
Thus we complete the proof of the proposition.
\end{proof}

The set $\{s_0, s_1, s_2\}$ form a~Coxeter group.
Let $G$ be the group generated by $s_0$, $s_1$ and~$s_2$.
The relations can be summarized as follows:

\begin{lemma}\label{lemma3}
The generators $s_0$, $s_1$, $s_2$ of the group $G$ satisfy the following relations:
\begin{gather}
\label{RelOurG1}
1)\quad
s_0^2 = s_1^2 = s_2^2 = \id,
\\
\label{RelOurG2}
2)\quad
(s_0 s_2)^2 = \id,
\qquad
(s_0 s_1)^4= (s_1 s_2)^4= \id,
\\
\label{RelOurG3}
3)\quad
(s_2 s_1 s_0 s_1)^3 = (s_1 s_2 s_1 s_0)^3 = \id.
\end{gather}
\end{lemma}

\begin{proof}
One can check by direct computation using the matrix realization given above.
We shall leave to readers.
\end{proof}

\begin{remark}
The relations~\eqref{RelOurG1} and~\eqref{RelOurG2} in Lemma~\ref{lemma3}
are the relations are same as that of the af\/f\/ine Weyl
group $W(\widetilde{C}_2)$:
\begin{center}
\setlength{\unitlength}{0.1mm}
\begin{picture}
(1200, 100) \put(600, 50){\circle{14}} \put(400, 50){\circle{14}} \put(800, 50){\circle{14}} \put(406,
46){\line(1,0){188}} \put(406, 54){\line(1,0){188}} \put(606, 46){\line(1,0){188}} \put(606, 54){\line(1,0){188}}
\put(390, 20){$s_0$} \put(590, 20){$s_1$} \put(790, 20){$s_2$}
\end{picture}

{\small Dynkin diagram of $\widetilde{C}_2$}
\end{center}
\end{remark}

We utilize properties of af\/f\/ine Weyl group $W(\widetilde{C}_2)$, especially translations in $W(\widetilde{C}_2)$ to
describe the structure of the group $G$ (For properties of af\/f\/ine Weyl groups, see Iwahori--Matsumoto~\cite{IM} and
Humphreys' book~\cite{HBook}).
Here we follow the notation of~\cite{HBook}.

In general, it is well known that af\/f\/ine Weyl group is a~semidirect product of a~Weyl group of the corresponding f\/inite
root system and the translation group corresponding to the coroot lattice.
We def\/ine the root vectors in the two dimensional Euclidean space~$V$ for the root system~$C_2$ as the following
picture:
\begin{center}
\setlength{\unitlength}{1mm}
\begin{picture}
(60, 60) \put(48, 27){{$ \alpha_1$}} \put(12, 49){{$ \alpha_2$}} \put(12, 17){{$ \alpha_0$}} \put(30, 30){\vector(1,
1){25}} \put(30, 30){\vector(0, 1){25}} \put(30, 30){\vector(-1, 0){25}} \put(30, 30){\vector(-1, -1){25}} \put(30,
30){\vector(0, -1){25}} \put(30, 30){\vector(1, -1){25}} \thicklines \put(30, 30){\vector(1, 0){25}} \put(30,
30){\vector(-1, 1){25}}
\end{picture}

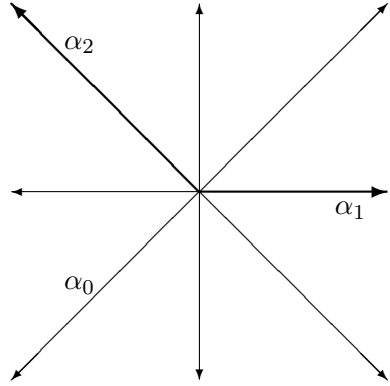
{\small Roots of root system $C_2$}

\end{center}
The null root $\alpha_0$ for $\widetilde{C}_2$ is given by $- 2\alpha_1-\alpha_2$.
For a~root $\alpha$, we denote by the corresponding coroot $\alpha^\vee$ given by ${\alpha^\vee = 2 \alpha
/ (\alpha, \alpha)}$, where $(\cdot, \cdot)$ is the Killing form.
In this case, the generators of the Weyl group of the root system $C_2$ be given by $s_1$ and $s_2$.
We denote $L^{\vee}$ by the coroot lattice of the root system $C_2$.
For $d \in V$, let $t(d)$ be the translation which sends $\lambda \in V$ to $d+\lambda$.

\begin{lemma}
The group $G$ is of order $72$.
Furthermore, $G$ is isomorphic to a~semidirect product of $W({C}_2)$ and $L^{\vee} /3L^{\vee}$.
\end{lemma}

\begin{proof}
Note that $s_2 s_1 s_0 s_1$ is the translation $t(\alpha_2^{\vee})$ of minimum length in $V$.
It is obvious to see that $s_1 t(\alpha_2^\vee) s_1 = s_1 s_2 s_1 s_0$ is $t(s_1 \alpha_2^\vee) = t(\alpha_1^\vee+\alpha_2^\vee)$.
Thus the group $G$ is a~subgroup of the group $W(C_2) \ltimes L^{\vee} / 3L^{\vee}$.
In order to see $G$ is isomorphic to $W(C_2) \ltimes L^{\vee} / 3L^{\vee}$, it suf\/f\/ices to check that $t(\alpha_2^\vee)
\not= \id$
and $t(\alpha_1^\vee+\alpha_2^\vee) \not= \id$.
Both of them can be done by direct computation using the matrix realization.
\end{proof}

\begin{remark}
The Coxeter group $G$ can be considered as a~``periodic'' af\/f\/ine Weyl group.
In particular, the relation~\eqref{RelOurG3} implies ``periodicity'' with respect to translations for the coroot lattice.
David Bessis informed us that the group~$G$ is {\it not} one of complex ref\/lection groups~\cite{ST}.
He proved this by calculating the character of the group~$G$.
\end{remark}

We are going to classify possible and non-trivial transformation for the $A_n$ ${}_4 F_3$ series of rectangular
type~\eqref{itSer}.

Let $H$ be the subgroup of the group $G$ generated by $s_0$ and $s_2$.
Recall that the ${}_4 F^n_3$ series of type~\eqref{itSer} is invariant under the action of $s_0$ and $s_2$.
Then our problem reduces to give an orbit decomposition of the double coset $H \backslash G / H$.
The representatives of this decomposition are given by
\begin{gather}
 (i)  \quad \id,
\qquad
 {(ii)} \quad s_1,
\qquad
 {(iii)} \quad s_1 s_2 s_1,
\qquad
 {(iv)} \quad s_1 s_0 s_1,
\qquad
 {(v)} \quad s_1 s_0 s_2 s_1,\nonumber
\\
{(vi)} \quad s_1 s_2 s_1 s_0 s_1,
\qquad
 {(vii)} \quad s_1 s_0 s_1 s_2 s_1,
\qquad
 {(viii)} \quad s_1 s_0 s_2 s_1 s_0 s_2 s_1.
\label{RepOurG}
\end{gather}

We are ready to exhibit a~complete list of possible transformations for the series of form~\eqref{itSer} according to
the representative~\eqref{RepOurG} of each orbit in~$G$.
We use Lemma~\ref{lemma2} frequently without stated otherwise in simplifying the factors.

The f\/irst one $(i)$ in~\eqref{RepOurG} is identical.
The second is~\eqref{ias2}.
The transformation correspon\-ding~$(iii)$ is
\begin{gather}
\SP{\{{-m_i}\}_n
\\
\{x_i\}_n}{a_1, a_2
\\
e_1, e_2}{c\\d} = \frac{[d+e_1 -a_2-c, e_1-a_1]_{|M|}}{[d+e-a_1-a_2-c, e_1]_{|M|}}
\nonumber
\\
\qquad{}
\times \SP{\{{-m_i}\}_n
\\
\{x_i\}_n}{a_1, d-c
\\d+e_2-a_2-c, d+e_1-a_2 -c}{d-a_2
\\d},
\label{ias1}
\end{gather}
which is equivalent to $q \to 1$ limit of the f\/irst $A_n$ Sears transformation~(4.23) of~\cite{KajiBDT}.
The one~$(iv)$ is
\begin{gather}
\SP{\{{-m_i}\}_n
\\
\{x_i\}_n}{a_1, a_2
\\
e_1, e_2}{c\\d}
= \frac{[d-c]_{|M|}}{[ d+e_1-a_1-a_2-c]_{|M|}} \prod\limits_{1 \le i \le n}\!\! \frac{[ d+e_1-a_1 -a_2+x_i ]_{m_i}}{[d+x_i ]_{m_i}}
\nonumber
\\
\qquad
{}
\times \SP{\{{-m_i}\}_n
\\
\{x_i\}_n}{e_1-a_1, e_1-a_2\\e_1+e_2-a_1-a_2, e_1}{c\\d+e_1-a_1-a_2}.
\label{ias3}
\end{gather}
The f\/ifth one $(v)$ is
\begin{gather}
\SP{\{{-m_i}\}_n
\\
\{x_i\}_n}{a_1, a_2
\\
e_1, e_2}{c\\d}
\nonumber\\
= \frac{[d+e_1-a_2-c, a_2]_{|M|}}{[ d+e-a_1 -a_2-c, e_1 ]_{|M|}} \prod\limits_{1 \le i \le n} \frac{[ d+e_1-a_1-a_2+x_i ]_{m_i}}{[d+x_i ]_{m_i}}
\label{iars2}
\\
\phantom{=}
\times
\SP{\{{-m_i}\}_n
\\
\{x_i\}_n}{e_1-a_2, d+e_1-a_1-a_2-c
\\
d +e_1+e_2-a_1-2 a_2- c, d+e_1-a_2-c}{d-a_2
\\
d+e_1-a_1-a_2}.
\nonumber
\end{gather}
The one $(vi)$ is
\begin{gather}
\SP{\{{-m_i}\}_n
\\
\{x_i\}_n}{a_1, a_2
\\
e_1, e_2}{c\\d}
\nonumber\\
= \frac{[ d-c, d+e_1+e_2-a_1-2 a_2 -c ]_{|M|}}{[ d+e_1-a_1-a_2-c, d+e_2-a_1-a_2 -c]_{|M|}}
\prod\limits_{1 \le i \le n} \frac{[ d-a_1+x_i]_{m_i}}{[ d+x_i ]_{m_i}}
\label{iars1}
\\
\times \SPa{\{{-m_i}\}_n
\\
\{\tilde{x}_i\}_n}{e_1-a_2, e_2-a_2
\\
e_1+e_2-a_1-a_2, d+e_1+e_2-a_1-2 a_2-c}
\!\!\!\!\left.
\left|
\begin{matrix}
e_1+e_2-a_1-a_2-c
\\
e_1+e_2-a_2 -c
\end{matrix}
\right|
1 \right).
\nonumber
\end{gather}
The seventh one $(vii)$ is
\begin{gather}
\SP{\{{-m_i}\}_n\\ \{x_i\}_n}{a_1, a_2\\e_1, e_2}{c\\d}
= \frac{[d+e_1-a_2-c, d+e_1-a_1-c]_{|M|}}{[ d+e-a_1-a_2-c, e_1 ]_{|M|}} \prod\limits_{1 \le i \le
n} \frac{[c+x_i ]_{m_i}}{[ d x_i ]_{m_i}}
\nonumber
\\
{}
\times
\SPa{\{{-m_i}\}_n
\\
\{\tilde{x}_i\}_n}{d -c, d+e_1-a_1-a_2-c\\d+e_1-a_1-c, d+e_1-a_2-c}
\!\!\!\left.
\left|
\begin{matrix}
e_1-c
\\
d+e_1+e_2-a_1-a_2-2 c
\end{matrix}
\right|
1 \right).\!\!\!
 \label{iars3}
\end{gather}
The one $(viii)$ is
\begin{gather}
\SP{\{{-m_i}\}_n
\\
\{{x}_i\}_n}{a_1, a_2
\\
e_1, e_2}{c\\d}
= \frac{[a_1, a_2]_{|M|}}{[e_1, d+e_1-a_1-a_2-c]_{|M|}} \prod\limits_{1 \le i \le n} \frac{[ c+x_i ]_{m_i}}
{[ d+x_i ]_{m_i}}
\nonumber\\
\qquad
\phantom{=}{}
\times \SPa{\{{-m_i}\}_n
\\
\{\tilde{x}_i\}_n}{d+e_1-a_1-a_2-c, d+e_2-a_1-a_2-c
\\
d +e_1+ e_2- a_1-2 a_2-c, d +e_1+ e_2-2 a_1- a_2-c}
\nonumber
\\
\hspace*{40mm}
\left.
\left|
\begin{matrix}
e_1+e_2-a_1-a_2-c
\\
d+e_1+e_2-a_1-a_2-2 c
\end{matrix}
\right|
1 \right).
\label{iara}
\end{gather}
\eqref{iara} has an alternative expression:
\begin{gather}
\SP{\{{-m_i}\}_n\\ \{x_i\}_n}{a_1, a_2\\e_1, e_2}{c\\d}
= \frac{[a_1, a_2]_{|M|}}{[e_1, e_2 ]_{|M|}} \prod\limits_{1 \le i \le n} \frac{[ c+x_i ]_{m_i}}{[ d+x_i]_{m_i}}
\nonumber
\\
\qquad
\times \SP{\{{-m_i}\}_n
\\
\{\tilde{x}_i\}_n}{1-|M|-e_1, 1-|M|-e_2
\\
1- |M| -a_1, 1-|M|-a_2}{1-|M|-d\\1-|M|-c}.
\label{iar}
\end{gather}
Note that this expression of the formula implies the reversing the order of the summation as~${}_4 F^n_3$ series of the
form~\eqref{itSer}.

\subsection[$A_n$ ${}_4 F_3$ series of triangular type]{$\boldsymbol{A_n}$ $\boldsymbol{{}_4 F_3}$ series of triangular type}

We describe the invariance group for triangular $A_n$ Whipple transformation~\eqref{las2}, namely the series of the form
\begin{gather*}
%\label{ltSer}
\SP{\{{b_i}\}_n
\\
\{x_i\}_n}{{-N}, a\\
e_1, e_2}{c\\d}
\end{gather*}
with the balancing condition
\begin{gather*}
%\label{lwbc}
a + B+c+{1-N} = d+e_1+e_2.
\end{gather*}

Note that, on contrast to the case of~\eqref{itSer}, the action of the permutation $s_0$
in Section~\ref{section2.3} is {\it not} valid.
So what we are to consider is the action of the permutation $s_2$ of the parame\-ters~$e_1$ and~$e_2$ and the
transformation of the parameters in~\eqref{las2}.
The action of each parameters of the transformation~\eqref{las2} is given by
\begin{gather*}
%\label{lMatrS}
s_1: \
\begin{bmatrix}
b
\\
a\\
c
\\
d
\\
e_1
\\
e_2
\end{bmatrix}
\to
\begin{bmatrix}
b
\\
e_1-a\\
e_1-c
\\
e_1+e_2-a~- c
\\
e_1
\\
d+e_1-a~- c
\end{bmatrix}
=S_1
\begin{bmatrix}
b
\\
a\\
c
\\
d
\\
e_1
\\
e_2
\end{bmatrix}
,
\qquad
S_1:=
\begin{bmatrix}
1 & 0 & 0 & 0 & 0 & 0
\\
0 & -1 & 0 & 0 & 1 & 0
\\
0 & 0 & -1 & 0 & 1 & 0
\\
0 & -1 & -1 & 0 & 1 & 1
\\
0 & 0& 0 & 0 & 1& 0
\\
0 & -1 & -1 & 1 & 1& 0
\end{bmatrix}
.
\end{gather*}

Let $G_t$ be the group generated by the transformations $s_1$ and $s_2$.
The relations between $s_1$ and $s_2$ are completely same as that between $s_1$ and $s_2$ in Section~\ref{section2.3}.
Namely,
\begin{gather*}
s_1^2 = s_2^2 = \id,
\qquad
(s_1 s_2)^4 = \id.
\end{gather*}
It follows that the group $G_t$ is isomorphic to $W(C_2)$, the Weyl group associated to the root system $C_2$.

To classify possible and non-trivial transformation formula, what is going to see is to give an orbit decomposition of the
double coset $H \backslash G_t /H$, where $H$ is a~subgroup of $G_t$ generated by $s_2$.
Note that $H$ is isomorphic to $\mathfrak{S}_2$.
The representatives of each orbits associated to this decomposition are given by $(i)$ $\id$, $(ii)$~$s_1$
and $(iii)$~$s_1 s_2 s_1$.

We now present the corresponding $A_n$ ${}_4 F_3$ transformations attached to each representative given above.
The transformation for $(i)$ is identical and the second one $(ii)$ is~\eqref{las2}.
The third one $(iii)$ is
\begin{gather}
\SP {\{b_i\}_n
\\
\{x_i\}_n}{{-N}, a\\
e_1, e_2}{c\\d} = \frac{[e_1-B, d+e_1-a~- c]_N}{[e, d+e_1-a~- B-c]_N}
\nonumber
\\
\qquad{}
\times \SP {\{b_i\}_n
\\
\{x_i\}_n}{{-N}, d- c
\\
d+e_2-a~- c, d+e_1-a~- c}{d-a~\\d},
\label{las1}
\end{gather}
which is the $q \to 1$ limit case of $A_n$ Sears transformation ((4.2) in~\cite{KajiS}).

\subsection{Remarks on the results of Section~\ref{section2}}
\label{section2.5}

Finally, we close the present paper to give remarks on the structure of the corresponding groups of the $A_n$
hypergeometric series of each cases.

\begin{remark}[the case when ${n=1}$ in ${A_n}$ ${{}_4 F_3}$ series]
In the case when $n=1$, all the transformations~\eqref{ias1},~\eqref{ias2} and~\eqref{ias3} of rectangular type
and~\eqref{las1} of triangular type reduce to the Whipple transformation formula~\eqref{d1st1}.
All the transformations~\eqref{iars1},~\eqref{iars2}, and~\eqref{iars3} of rectangular type reduce to~\eqref{d1rst1}.
The transformation~\eqref{iar} of rectangular type reduces to~\eqref{ad1r1} and implies the reversing of the order of
the summation in the $A_n$ ${}_4 F_3$ series of rectangular type.
\end{remark}

\begin{remark}[correspondence of the groups in Sections~\ref{section2.2} and~\ref{section2.3}]
By direct manipulation of the matrix realization in Section~\ref{section2.3},
we have $s=\sigma_4 \sigma_3 \sigma_1 \sigma_5 \sigma_4$.
Thus we f\/ind that the group $G$ is isomorphic to the subgroup of~$\mathfrak{S}_6$ generated by~$\sigma_2$,~$\sigma_5$ and~$s$.
\end{remark}

\begin{remark}
Except for Hardy type invariant form~\eqref{R4F3} for ${}_4 F_3^n$ series of the form~\eqref{itSer}, all other results
are valid in the basic case and one can obtain in the same line as the discussion in this section.
For Hardy type invariant form for terminating balanced~${}_4 \phi_3$ series have already appeared in Van der Jeugt and
Srinivasa Rao~\cite{vdj}.
\end{remark}

\section[Symmetry groups of $A_n$ elliptic hypergeometric series]{Symmetry groups
of $\boldsymbol{A_n}$ elliptic hypergeometric series}
\label{section3}

\subsection[Preliminaries on $A_n$ elliptic hypergeometric series]{Preliminaries
on $\boldsymbol{A_n}$ elliptic hypergeometric series}

Here, we give notations for (multiple) elliptic hypergeometric series and recall the results of our previous paper with
M.~Noumi~\cite{KajiNou}.

Let $[[x]]$ be a~non-zero and homomorphic odd function in $\mathbb C$ which satisf\/ies the Riemann relation:
\begin{gather}
1)
\quad
[[-x]]=-[[x]],
\nonumber
\\
2)
\quad
[[x+y]]\,[[x-y]]\,[[u+v]]\,[[u-v]]
\nonumber
\\
\phantom{2)}
\quad
\qquad
=[[x+u]]\,[[x-u]]\,[[y+v]]\,[[y-v]]
-[[x+v]]\,[[x-v]]\,[[y+u]]\,[[y-u]].
\label{RieRel}
\end{gather}
There are following three classes of such functions:
\begin{itemize}\itemsep=0pt
\item
$\sigma(x;\omega_1,\omega_2)$: Weierstrass sigma function with the periods $(\omega_1,\omega_2)$ (elliptic),
\item
$\sin(\pi x)$: the sine function (trigonometric),
\item
$x$: rational.
\end{itemize}
It is classically known~\cite{WW} that all function $[[x]]$ satisfy the condition~\eqref{RieRel} are obtained from above
three functions by transformation of the form $e^{ax^2+b} [[ c x ]]$ for complex numbers $a, b, c \in {\mathbb C}$.

Fix a~generic constant $\delta\in\mathbb{C}$ so that for all integer $k\in\mathbb{Z}$, $[[k\delta]]$ does not equal to
zero.
In the case when $[[x]]$ is Weierstrass sigma function $\sigma(x;\omega_1,\omega_2)$ (the elliptic case for short),
the condition for $\delta$ is given by $\delta\not\in\mathbb{Q}\omega_1+\mathbb{Q}\omega_2$.

Throughout the present paper, we consider the function $[[ x ]]$ as the elliptic case unless other\-wise stated.

Next a~shifted factorial $[[x]]_k$ is def\/ined by
\begin{gather*}
[[x]]_k:=[[x]][[x+\delta]]\cdots [[x+(k-1)\delta]],
\qquad
 k=0,1,2,\ldots.
\end{gather*}
Further, we denote
\begin{gather*}
[[x_1, \dots, x_r]]_k:= [[x_1]]_k \cdots [[x_r]]_k.
\end{gather*}

The elliptic hypergeometric series ${}_{r+3}E_{r+2}$ is def\/ined as follows
\begin{gather*}
%\label{DefEHS}
{}_{r+3}E_{r+2}(s; \{u_k\}_{r}) = {}_{r+3}E_{r+2}(s; u_1, \dots,u_r)
:=
\sum\limits_{m \in {\mathbb N}}
\frac{[[s+2 m \delta]]}{[[s]]} \frac{[[s]]_m}{[[\delta]]_m}
\prod\limits_{1 \le i \le r}\frac{[[u_i]]_m}{[[\delta+s-u_i]]_m}.
\end{gather*}
In the case when $[[x]]$ is a~trigonometric function $\sin x$,
this series reduces to the basic very well-poised hypergeometric series ${}_{r+3} W_{r+2}$.
Note that ${}_{r+3} E_{r+2}$ series are also symmetric with respect to the parameter $u_k$ for $1 \le k \le r$.

All the ${}_{r+3} E_{r+2}$ series discussed in this paper is balanced, namely we assume
\begin{gather*}
%\label{BalanceE}
u_1+\dots+u_{r} = {\frac{r-1}{2}} s+\frac{r-3}{2}.
\end{gather*}

Now, we note the conventions for naming series as~$A_n$ elliptic hypergeometric series (or referred as elliptic
hypergeometric series in ${\rm SU}(n+1)$).
Let $\gamma= (\gamma_1, \dots, \gamma_n) \in \mathbb{N}^n$ be a~multi-index.
We denote generalizations of the Vandermonde determinant
\begin{gather*}
\Delta[x]:= \prod\limits_{1 \le i < j \le n} [[x_i-x_j]]
\qquad
\mbox{and}
\qquad
\Delta[x+\gamma \delta]:= \prod\limits_{1 \le i < j \le n} [[x_i+\gamma_i \delta-x_j-\gamma_j \delta ]],
\end{gather*}
for the sets of variables $x=(x_1, \dots, x_n)$ and $x+\gamma \delta = (x_1+{\gamma_1} \delta, \dots, x_n
+{\gamma_n} \delta)$ respectively.
In this paper we refer multiple series of the form
\begin{gather}
\label{DefAnEHS}
\sum\limits_{\gamma \in {\mathbb N}
^n} \frac{\Delta [x +{\gamma \delta}]}{\Delta [x]} H(\gamma)
\end{gather}
which reduce to elliptic hypergeometric series ${}_{r+1} E_r$ for a~nonnegative integer $r$ when $n=1$ and symmetric
with respect to the subscript $1 \le i \le n$ as~$A_n$ elliptic hypergeometric series.
Other terminology are similar to the case of $A_n$ (ordinary) hypergeometric series.
The subscript $n$ in the label $A_n$ attached to the series is the dimension of the multiple series~\eqref{DefAnEHS}.

We are going to introduce the multiple elliptic hypergeometric series $E^{n,m}$ which is de\-f\/i\-ned~by
\begin{gather*}
%\label{MEWP}
\ME{n, m}{\{a_i\}_n
\\
\{x_i\}_n}{s}{\{u_k\}_{m}}{\{v_k\}_{m}}{}
:=
\sum\limits_{\gamma \in \mathbb{N}
^n}{} \frac{\Delta[x+\gamma \delta]}{\Delta[x]}\prod\limits_{1 \le i \le n} \frac{[[({|\gamma|+\gamma_i}) \delta +s
+ x_i]]}{[[s+x_i]]}
\\
\qquad
\phantom{:=}
\times
\prod\limits_{1\le j \le n} \frac{[[s+x_j]]_{|\gamma|}}{[[\delta+s-a_j+x_j ]]_{|\gamma|}}
\left(\prod\limits_{1 \le i \le n} \frac{[[a_j+x_i-x_j]]_{\gamma_i}}{[[ \delta+x_i-x_j ]]_{\gamma_i}} \right)
\\
\qquad
\phantom{:=}
\times
\prod\limits_{1 \le k \le m} \frac{[[v_k]]_{|\gamma|}}{[[ \delta+s-u_k]]_{|\gamma|}} \left(\prod\limits_{1
\le i \le n} \frac{[[u_k+x_i]]_{\gamma_i}}{[[ \delta+s-v_k+x_i]]_{\gamma_i}} \right).
\end{gather*}

In the case when $n=1$, $E^{1, m}$ series reduces to (one dimensional) elliptic hypergeometric series ${}_{2m+4} E
_{2m+3} \left(s; \{u_k\}_m, \{v_k\}_m \right)$.

Note that $E^{n,m}$ series is symmetric within two sets of parameters $\{u_k\}_m$ and $\{v_k\}_m$ respectively.
This fact will be a~key of the latter discussion of the symmetry for the $E^{n,m}$ series.

Here, we present the balanced duality transformation formula for multiple elliptic hypergeometric series
from~\cite{KajiNou}.

Under the balancing condition
\begin{gather*}
%\label{bcEBDT}
c_1+c_2+d_1+\sum\limits_{1 \le i \le n}
a_i+\sum\limits_{1 \le k \le m}
(u_k+v_k) = (m+N +1) \delta+(m+2) s.
\end{gather*}
We have the following transformation formula between $E^{n, m+2}$ ($A_n$ ${}_{2m+ 8} E_{2m+7}$) series and $E^{m,
n+2}$ ($A_m$ ${}_{2n+8} E_{2n+7}$) series ((3.17) in~\cite{KajiNou}):
\begin{gather}
\ME{n,m+2}{\{a_i\}_n
\nonumber
\\
\{x_i\}_n}{s}{c_1,c_2,\{u_k\}_m}{d_1, -N \delta, \{v_k\}_m}
\\
\qquad
= \frac{[[\delta+s-c_1-d_1, \delta+s-c_2-d_1]]_N}{[[\delta+s-c_1, \delta+s-c_2]]_N} \prod\limits_{1 \le k \le m}
\frac{[[v_k, \delta+s-u_k-d_1]]_N}{[[\delta+s-u_k, v_k-d_1]]_N}
\nonumber
\\
\qquad
\phantom{=}
\times \prod\limits_{1 \le i \le n} \frac{[[\delta+s+x_i, \delta+s+x_i-a_i-d_1]]_N}{[[\delta+s+x_i-a_i,
\delta+s+x_i-d_1]]_N}
\nonumber
\\
\qquad
\phantom{=}
\times \ME{m,n+2}{\{b_k\}_m
\\
\{y_k\}_m}{t}{-c_1,-c_2, \{z_i\}_n}{d_1,- N \delta, \{w_i\}_n},
\label{EBDT1}
\end{gather}
where
\begin{gather*}
t=d_1+d_2-s-\delta,
\qquad
b_k=\delta+s-u_k-v_k,
\qquad
 y_k=\delta+s-v_k,
 \qquad
 k=1,\ldots,m,
\\
z_i=x_i-a_i,
\qquad
 w_i=d_1+d_2-s-x_i,
 \qquad
 i=1,\ldots,n.
\end{gather*}

\begin{remark}
In~\cite{RoseKaji}, Rosengren also obtained~\eqref{EBDT1} by a~dif\/ferent way from~\cite{KajiNou}.
In the case when $m=n=1$ and $x_1 = y_1 =0 $,~\eqref{EBDT1} reduces to the following transformation formula for
terminating balanced ${}_{10} E_9$ series:
\begin{gather}
{}_{10} E_{9} \left(s; c_0,c_1,c_2,c_3,d_0,d_1, -N \delta \right) = \frac{[[d_0, \delta+s]]_N}{[[d_0-d_1,
\delta+s-d_1]]_N}
\nonumber
\\
\qquad{}
\times \prod\limits_{0 \le k \le 3} \frac{[[\delta+s-c_k-d_1]]_N}{[[\delta+s-c_k]]_N}{}_{10} E_9
\big(
\widetilde{s}; \widetilde{c}_0,\widetilde{c}_1, \widetilde{c}_2, \widetilde{c}_3,\widetilde{d}_0,d_1,-N \delta
\big),
\label{mn1EBDT1}
\\
(c_0+c_1+c_2+c_3+d_0+d_1=(2+N) \delta+3 s),
\nonumber
\end{gather}
where
\begin{gather*}
\widetilde{s}=d_1+d_2-d_0,
\qquad
\widetilde{d}_0=d_1+d_2-s,
\qquad
\widetilde{c}_k=\delta+s-d_0-c_k,
\qquad
 k=0,1,2,3.
\end{gather*}

Note that the ${}_{10} E_9$ transformation~\eqref{mn1EBDT1} can also be obtained by iterating twice in an appropriate
manner the (rather well-known) elliptic version of the Bailey transformation~\cite{B1} (see also~\cite{BB}) for ${}_{10} E_9$ series due to Frenkel and
Turaev~\cite{FT}:
\begin{gather}
{}_{10}E_{9} \left({s}; {c}_0,{c}_1,{c}_2, d_0,d_1,d_2, -N \delta \right) = \frac
{[[\delta+s]]_N}{[[\delta+s-d_0-d_1-d_2]]_N}\nonumber
\\
\qquad{}
\times \prod\limits_{0 \le k \le 2} \frac {[[\delta+s-d_0-d_1-d_2+ d_k]]_N}{[[\delta+s-d_k]]_N}{}_{10}E_{9}
\left(
\widetilde{s}; \widetilde{c_0}, \widetilde{c_1}, \widetilde{c_2}, d_0,d_1,d_2,-N \delta \right),
\label{EBaileyT1}
\\
\widetilde{s}=\delta+2 s-c_0-c_1-c_2,
\qquad
\widetilde{c}_0=\delta+s-c_1-c_2,
\qquad
\widetilde{c}_1=\delta+s-c_0-c_2,
\nonumber\\
\widetilde{c}_2=\delta+s-c_0-c_1,
\qquad
c_0+c_1+c_2+d_0+d_1+d_2=(2+N) \delta+3 s.
\nonumber
\end{gather}
Similarly,~\eqref{EBaileyT1} can also be obtained by iterating~\eqref{mn1EBDT1} (see~\cite{KajiNou} and latter
discussions).

In the case when $m = 1$, $y_1 = 0 $,~\eqref{EBDT1} reduces to the transformation formula between \mbox{$n$-dimensional}
$E^{n,3}$ series and $1$-dimensional ${}_{2n+8}E_{2n+7}$ series
\begin{gather}
\ME{n,3}{\{a_i\}_n
\\
\{x_i\}_n}{s}{c_0,c_1,c_2}{d_0,d_1, -N \delta}
\nonumber \\
\qquad{}
= \frac{[[d_0, \delta+s-c_0-d_1, \delta+s-c_1-d_1, \delta+s-c_2-d_1]]_N}{[[d_0-d_1, \delta+s-c_0, \delta+s-c_1,
\delta+s-c_2]]_N}
\nonumber
\\
\qquad
\phantom{=}{}
\times \prod\limits_{1 \le i \le n} \frac{[[\delta+s+x_i, \delta+s+x_i-a_i-d_1]]_N}{[[\delta+s+x_i-a_i,
\delta+s+x_i-d_1]]_N}
\nonumber
\\
\qquad
\phantom{=}{}
\times
{}_{2 n+8}E_{2 n+7} \left(t; e_0, e_1,e_2, \{u_i\}_n, \{v_i\}_n, d_1,-N \delta \right),
\label{m1EBDT}
\end{gather}
where
\begin{gather*}
t=d_1-N \delta -d_0,
\qquad
 e_k=\delta+s-d_0-c_k,
\qquad
 k=0,1,2 ,
\\
u_i=\delta+s-d_0+x_i-a_i,
\qquad
 v_i=d_1-N \delta -s-x_i,
 \qquad
 i=1,\ldots,m ,
\end{gather*}
under the balancing condition for $E^{n,3}$ series
\begin{gather*}
\sum\limits_{1 \le i \le m}a_i+ (c_0+c_1+c_2)+(d_0+d_1) = (2+N) \delta+3s.
\end{gather*}
\end{remark}

\subsection[Symmetry of ${}_{10} E_9$ series ($A_1$ case)]{Symmetry
of $\boldsymbol{{}_{10} E_9}$ series ($\boldsymbol{A_1}$ case)}
\label{section3.2}

Here, we describe the symmetry of $1$-dimensional elliptic Bailey transformation for ${}_{10} E_9$
series~\eqref{EBaileyT1}, namely for the ${}_{10} E_9$ series of the form
\begin{gather*}
%\label{10E9}
{}_{10}E_{9} \left({s}; {c}_0,{c}_1,{c}_2, c_3,c_4, c_5, -N \delta \right),
\end{gather*}
with the balancing condition
\begin{gather}
\label{1bc}
c_0+c_1+c_2+c_3+c_4+c_5= (2+N) \delta+3s.
\end{gather}
The result here has appeared in the paper by S.~Lievens and J.~Van der Jeugt~\cite{LvdJ} in the case of very well-poised
basic hypergeometric series ${}_{10} W_9$.
But our description given here is modif\/ied for the sake of the connection of the results in this section.

For $k=1, \dots, 5$, let $s_k$ be the permutation for the parameters $c_{k-1}$ and $c_k$.
Let $b$ be the transformation of parameters for the Bailey transformation~\eqref{EBaileyT1}.
Note that these are af\/f\/ine transformations in $7$-dimensional vector space.
Here we shall give a~$8 \times 8$ matrix realization acting on the vector $\vec{v}_1 = {}^t[s, c_0, c_1, c_2, c_3, c_4,
c_5, \delta]$ for these transformations.
The transformation of parameters $b$ for Bailey transformation~\eqref{EBaileyT1} and its matrix realization are given by
\begin{gather*}
b: \ \vec{v}_1 =
\begin{bmatrix}
s^{}
\\
c_0^{}
\\
c_1^{}
\\
c_2^{}
\\
c_3^{}
\\
c_4^{}
\\
c_5^{}
\\
\delta
\end{bmatrix}
\mapsto
\begin{bmatrix}
2s+\delta-c_0 -c_1 -c_2
\\
s+\delta-c_1-c_2
\\
s+\delta-c_0-c_2
\\
s+\delta-c_0-c_1
\\
c_3
\\
c_4
\\
c_5
\\
\delta
\end{bmatrix}
= B \cdot \vec{v}_1,
\end{gather*}
and
\begin{gather*}
B=
\begin{bmatrix}
2 & -1 & -1 & -1 & 0 & 0 & 0 & 1
\\
1 & 0 & -1 & -1 & 0 & 0 & 0 & 1
\\
1 & -1 & 0 & -1 & 0 & 0 & 0 & 1
\\
1 & -1 & -1 & 0 & 0 & 0 & 0 & 1
\\
0 & 0 & 0 & 0 & 1 & 0 & 0 & 0
\\
0 & 0 & 0 & 0 & 0 & 1 & 0 & 0
\\
0 & 0 & 0 & 0 & 0 & 0 & 1 & 0
\\
0 & 0 & 0 & 0 & 0 & 0 & 0 & 1
\end{bmatrix}
,
\end{gather*}
respectively.
The matrix realization for $s_1$ is given by
\begin{gather*}
S_1=
\begin{bmatrix}
1 & 0 & 0 & 0 & 0 & 0 & 0 & 0
\\
0 & 0 & 1 & 0 & 0 & 0 & 0 & 0
\\
0 & 1 & 0 & 0 & 0 & 0 & 0 & 0
\\
0 & 0 & 0 & 1 & 0 & 0 & 0 & 0
\\
0 & 0 & 0 & 0 & 1 & 0 & 0 & 0
\\
0 & 0 & 0 & 0 & 0 & 1 & 0 & 0
\\
0 & 0 & 0 & 0 & 0 & 0 & 1 & 0
\\
0 & 0 & 0 & 0 & 0 & 0 & 0 & 1
\end{bmatrix}
,
\end{gather*}
and so on.

\begin{proposition}[Hardy type invariant form for ${{}_{10} E_9}$ series with the condition~\eqref{1bc}]
If~the balancing condition~\eqref{1bc} holds,
\begin{gather*}
{}_{10} \widetilde{E}_{9} \left((\vec{v}_1)\right) :=
\frac{\prod\limits_{0 \le k \le 5}[[ \delta+s-c_k ]]_N}
{[[s]]_N}{}_{10}E_{9} \left({s}; {c}_0,{c}_1,{c}_2, c_3,c_4, c_5, -N \delta \right)
\end{gather*}
is invariant under the action of $b$ and $s_k$ for all $1\le k \le 5$.
\end{proposition}

Note that $B^2 = \id_8$, namely $b^2 = \id$.
By def\/inition, $\{s_i \,|\, i=1, 2, 3, 4, 5\}$ satisfy the relation
\begin{gather}
\label{rels}
\begin{cases}
s_i \not= \id,
\qquad
s_i^2 = \id, &  i = 1, 2, 3, 4, 5 ,
\\
s_{i} s_{i+1} s_{i} = s_{i+1} s_{i} s_{i+1}, &  i = 1, 2, 3, 4 ,
\\
s_{i}s_{j} = s_{j} s_{i}, &  |i-j| \ge 2 .
\end{cases}
\end{gather}
The relations between $b$ and $s_i$, $i=1, 2, 3, 4, 5$ summarized as the following lemma:

\begin{lemma}
We have
$
(s_3 b)^3 = \id$.
And $b$ commutes with any other $s_i$, $i=1, 2, 4, 5$.
\end{lemma}

The proof of this lemma can be given by using the matrix realization given above.

By combining this lemma and~\eqref{rels}, we conclude:

\begin{theorem}
The group $G_1$ generated by the permutations of the parameters $\{c_k \,| \,k=0, \dots, 5\}$ and the Bailey
transformation~\eqref{EBaileyT1} for ${}_{10} E_9$ series is isomorphic to the Weyl group $W(E_6)$ associated to the
root system~$E_6$:
\begin{center}
\setlength{\unitlength}{0.1mm}
\begin{picture}
(1200, 300) \put(600, 50){\circle{14}} \put(200, 50){\circle{14}} \put(400, 50){\circle{14}} \put(800, 50){\circle{14}}
\put(1000, 50){\circle{14}} \put(600, 250){\circle{14}} \put(207, 50){\line(1,0){186}} \put(407, 50){\line(1,0){186}}
\put(607, 50){\line(1,0){186}} \put(807, 50){\line(1,0){186}} \put(600, 57){\line(0,1){186}} \put(190, 20){$w_1$}
\put(390, 20){$w_3$} \put(590, 20){$w_4$} \put(790, 20){$w_5$} \put(990, 20){$w_6$} \put(610, 240){$w_2$}
\end{picture}

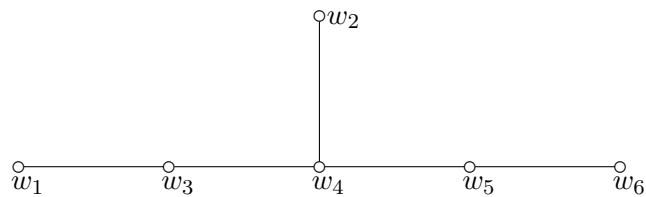
{\small Dynkin diagram of $E_6$}
\end{center}
\end{theorem}

Here we classify the possible non-trivial transformation formulas for ${}_{10} E_9$ series.
Notice again that ${}_{10} E_9$ series is symmetric for the permutation of parameters $c_0, \dots, c_5$.
That is, it is invariant under the action $s_i$ for all $i = 1, \dots, 5$.
Thus our problem turn out to give an orbit decomposition of the double coset $\mathfrak{S}_6 \backslash W (E_6) /
\mathfrak{S}_6$.

We def\/ine the mapping $\pi_1$ according to Bourbaki~\cite{Bour} as follows
\begin{gather*}
s_1 \mapsto w_1, \qquad b \mapsto w_2, \qquad s_i \mapsto w_{i+1}, \qquad  i=2, 3, 4, 5, 6 .
\end{gather*}

The representatives of orbits in the double coset $\mathfrak{S}_6 \backslash W (E_6) / \mathfrak{S}_6$ are given as
follows:
\begin{gather*}
1) \quad \tau_1 = \id,
\\
2) \quad \tau_2= w_2,
\\
3) \quad \tau_3 = w_2 w_4 w_3 w_5 w_4 w_2,
\\
4) \quad \tau_4 = w_2 w_4 w_3 w_1 w_5 w_4 w_3 w_6 w_5 w_4 w_2,
\\
5) \quad\tau_5= w_2 w_4 w_3 w_1 w_5 w_4 w_2 w_3 w_4 w_5 w_6 w_5 w_4 w_2 w_3 w_1 w_4 w_3 w_5 w_4 w_2.
\end{gather*}

Thus we are ready and we shall exhibit a~list of the possible ${}_{10} E_9$ transformations.
We assume that all the ${}_{10} E_9$ series in the formulas listed here satisfy the balancing condition~\eqref{1bc}.

The transformation corresponding to $\tau_1 = \id$
is identical as~${}_{10} E_9$ transformation.
The transformation corresponding to $\tau_2$ is equivalent to the Bailey transformation due to
Frenkel--Turaev~\eqref{EBaileyT1}.
The third one corresponding to $\tau_3$ is equivalent to~\eqref{mn1EBDT1}.
The forth one ($\tau_4$)~is
\begin{gather}
{}_{10}E_{9} \left({s}; {c}_0, {c}_1, {c}_2, c_3,c_4,c_5, -N \delta \right)
= \frac{[[\delta+s]]_{N}}{[[3\delta+4 s -c_0 -c_1-c_2-2 c_3-2 c_4-2 c_5 ]]_{N}}
\nonumber
\\
\qquad{}
\times \prod\limits_{0 \le k \le 2} \frac {[[c_{k+3}, \delta+s-c_0-c_1-c_2+c_k ]]_N}{[[\delta+s-c_k,
\delta+s-c_{k+3}]]_{N}}
\,  {}_{10}E_{9} \left(\widehat{s}; \widehat{c}_0, \widehat{c}_1,\widehat{c}_2, \widehat{d}_0, \widehat{d}_1,
\widehat{d}_2, -N \delta \right),\label{BaileyT3}
\\
\widehat{s}= (1-N) \delta+ s-c_3-c_4-c_5 = 3\delta+ 4 s- c_0-c_1-c_2-2 c_3 -2 c_4- 2 c_5,
\nonumber
\\
\widehat{c}_k=\delta+s-c_3 -c_4 -c_5+c_{k+3},
\nonumber
\\
\widehat{c}_{k+3} = -N \delta+c_k-s = 2\delta+ 2 s- c_0-c_1-c_2 -c_3-c_4-c_5+c_k,
\qquad
 k=0, 1, 2 .\nonumber
\end{gather}

Finally, the f\/ifth one corresponding to $\tau_5$ is
\begin{gather}
{}_{10} E_{9} \left({s}; {c}_0,{c}_1,{c}_2, c_3,
c_4, c_5, -N \delta \right)
= \frac{[[\delta+s]]_{N}}{[[ 4\delta+5 s-2 c_0- 2 c_1- 2 c_2- 2 c_3- 2 c_4 -2 c_5]]_{N}}
\nonumber
\\
\qquad
\phantom{=}{}
\times \prod\limits_{0 \le k \le 5} \frac{[[c_k]]_{N}}{[[\delta+s-c_k]]_{N}}{}_{10} E_{9} \left(\check{s};
\check{c}_0,\check{c}_1,\check{c}_2, \check{c}_3,\check{c}_4,\check{c}_5, -N \delta \right),
\label{BaileyT4}
\\
\check{s}= -2 N \delta-s =4\delta+5 s-2 c_0-2 c_1-2 c_2-2 c_3-2 c_4-2 c_5,
\nonumber
\\
\check{c}_k=-N \delta+c_k-s = 2\delta+2 s-c_0-c_1-c_2-c_3-c_4-c_5+c_k,
\qquad
 k= 0, 1, 2, 3, 4, 5 .
\nonumber
\end{gather}

Note that~\eqref{BaileyT4} implies reversing order of the summation in ${}_{10} E_9$ series.

\subsection[Symmetry of $A_n$ Bailey transformations of rectangular type]{Symmetry
of $\boldsymbol{A_n}$ Bailey transformations of rectangular type}
\label{section3.3}

{\sloppy Here we discuss the symmetry for two $A_n$ elliptic Bailey transformation formulas~\eqref{MN} and~\eqref{KN}.
The corresponding series is $A_n$ elliptic hypergeometric series of rectangular type, which the multiple series
terminates with respect to a~multi-index.
Namely the $E^{n, 3}$ series of the form
\begin{gather}
\label{MTEn}
\ME{n,3}{\{- m_i \delta\}_n
\\
\{x_i\}_n}{s}{c_0, c_1, c_2}{d_0,d_1,d_2}
\end{gather}
with the balancing condition
\begin{gather}
\label{bc}
(c_0+c_1+c_2)+(d_0+d_1+d_2) = (2+|M|) \delta+3s.
\end{gather}

}

In~\cite{KajiNou}, we obtained several $A_n$ generalizations of the elliptic Bailey transformation
formula~\eqref{EBaileyT1} for $E^{n, 3}$ series by iterating~\eqref{m1EBDT} twice.
Among these, here we give two transformations which $E^{n, 3}$ series of rectangular type which we discuss here.
These can be obtained in a~similar way as in Section~\ref{section3.2}.

{\bf $\boldsymbol{A_n}$ Bailey transformation for $\boldsymbol{E^{n, 3}}$ of rectangular type~(\ref{MTEn})} (Theorem~4.2
in~\cite{KajiNou}).
Suppose that $a_i =-m_i \delta$, $m_i \in {\mathbb N}$ for all $i=1, \dots, n$.
For $c_k$, $d_k$, $k=0,1,2$, suppose that the balancing condition~\eqref{bc}.
Then we have two types of $A_n$ Bailey transformation formula.

  {\bf $\boldsymbol{A_n}$ Bailey~I (Milne--Newcomb type)}
\begin{gather}
\ME{n,3}{\{- m_i \delta\}_n
\\
\{x_i\}_n}{s}{c_0,c_1,c_2}{d_0,d_1,d_2}
= \frac {[[\delta+s-c_1-d_0, \delta+s-c_2-d_0]]_{|M|}}{[[\delta+s-c_1, \delta+s-c_2]]_{|M|}}
\nonumber
\\
\qquad
\phantom{=}{}
\times \prod\limits_{1 \le i \le n} \frac{[[\delta+s+x_i, 2\delta +2 s -c_0-d_0-d_1-d_2 +x_i]]_{m_i}}
{[[\delta+s-d_0 +x_i, 2\delta +2 s-c_0-d_1-d_2 +x_i]_{m_i}}
\nonumber
\\
\qquad
\phantom{=}{}
\times \ME{n,3}{\{- m_i \delta\}_n
\\
\{x_i\}_n}{\widetilde{s}}{\widetilde{c}_0,c_1,c_2}{d_0, \widetilde{d}_1,\widetilde{d}_2},
 \label{MN}
\end{gather}
where
\begin{gather*}
%\label{MNP}
\widetilde{s}=\delta+2 s-c_0-d_1-d_2,
\qquad
\widetilde{c}_0=\delta+s-d_1-d_2,
\\
\widetilde{d}_1=\delta+s-c_0-d_2,
\qquad
\widetilde{d}_2=\delta+s-c_0-d_1.
\end{gather*}

  {\bf $\boldsymbol{A_n}$ Bailey~II (Kajihara--Noumi type)}
\begin{gather}
\ME{n,3}{\{- m_i \delta\}_n
\\
\{x_i\}_n}{s}{c_0, c_1, c_2}{d_0,d_1,d_2}
= \prod\limits_{1 \le i \le n} \left[ \frac{[[\delta+s+x_i, \delta+s -d_0-d_1 +x_i]]_{m_i}}{[[\delta+s-d_0 +x_i,
\delta+s-d_1 +x_i]]_{m_i}} \right.
\nonumber
\\
\qquad
\phantom{=}{}
\times \left.
\frac{[[ \delta+s-d_0-d_2 +x_i, \delta+s-d_1-d_2+x_i]]_{m_i}}{[[\delta+s-d_2 +x_i, \delta+s-d_0-d_1-d_2
+x_i]]_{m_i}} \right]
\nonumber
\\
\qquad
\phantom{=}{}
\times \ME{n,3}{\{- m_i \delta\}
\\
\{\widetilde{x}_i\}}{\widetilde{s}}{\widetilde{c}_0,\widetilde{c}_1, \widetilde{c}_2}{d_0,d_1,d_2},
\label{KN}
\end{gather}
where
\begin{gather*}
%\label{KNP}
\widetilde{s}=\delta+2 s-c_0-c_1-c_2,
\qquad
\widetilde{c}_0=\delta+s-c_1-c_2,
\qquad
\widetilde{c}_1=\delta+s-c_0-c_2,
\\
 \widetilde{c}_2=\delta+s-c_0-c_1,
\qquad
\widetilde{x}_i=- m_i \delta -x_i+|M| \delta,
\qquad
 i=1,\ldots,n.
\end{gather*}

\begin{remark}
In the case when $n=1$, $x_1 = 0$, both~\eqref{MN} and~\eqref{KN} reduce to the elliptic Bailey transformation for
${}_{10} E_9$ series~\eqref{EBaileyT1}.
Bailey~I~\eqref{MN} is originally due to Rosengren~\cite{RoseE} which is a~elliptic version of $A_n$ Bailey
transformation formula by Milne--Newcomb~\cite{MilNew1}.
Bailey~II~\eqref{KN} has originally appeared in our previous work~\cite{KajiNou} together with the basic case.
\end{remark}

First, we shall show the invariance property for the transformations for $E^{n,3}$ series of rectangular
type~\eqref{MTEn}.
Suppose that $n \ge 2$ till we will state otherwise.
Recall that the transformations of coordinates in the right hand side of the Bailey~I~\eqref{MN} and Bailey~II~\eqref{KN} are described as follows
\begin{gather*}
b_1: \
\begin{bmatrix}
s^{}
\\
c_0^{}
\\
c_1^{}
\\
c_2^{}
\\
d_0^{}
\\
d_1^{}
\\
d_2^{}
\end{bmatrix}
\mapsto
\begin{bmatrix}
s^{1, 0}
\\
c_0^{1, 0}
\\
c_1^{1, 0}
\\
c_2^{1, 0}
\\
d_0^{1, 0}
\\
d_1^{1, 0}
\\
d_2^{1, 0}
\end{bmatrix}
=
\begin{bmatrix}
2 s+\delta-c_0-d_1-d_2
\\
s+\delta-d_1-d_2
\\
c_1
\\
c_2
\\
d_0
\\
s+\delta-c_0-d_2
\\
s+\delta-c_0-d_1
\end{bmatrix}
,
\\
b_2: \
\begin{bmatrix}
s^{}
\\
c_0^{}
\\
c_1^{}
\\
c_2^{}
\\
d_0^{}
\\
d_1^{}
\\
d_2^{}
\end{bmatrix}
\mapsto
\begin{bmatrix}
s^{0, 1}
\\
c_0^{0, 1}
\\
c_1^{0, 1}
\\
c_2^{0, 1}
\\
d_0^{0, 1}
\\
d_1^{0, 1}
\\
d_2^{0, 1}
\end{bmatrix}
=
\begin{bmatrix}
2s+\delta-c_0 -c_1 -c_2
\\
s+\delta-c_1-c_2
\\
s+\delta-c_0-c_2
\\
s+\delta-c_0-c_1
\\
d_0
\\
d_1
\\
d_2
\end{bmatrix}
.
\end{gather*}

Note that these are compositions of linear transformations for parameters $s$, $c_0$, $c_1$, $c_2$, $d_0$, $d_1$, $d_2$ and shift
by $\delta$, namely af\/f\/ine transformations of 7-dimensional vector space.
Thus we give a~realization for these transformations in terms of $8 \times 8$ matrices acting on the vector $\vec{v} =
{}^t [s,c_0, c_1, c_2, d_0, d_1, d_2, \delta ]$ as follows
\begin{gather*}
\begin{bmatrix}
s^{1, 0}
\\
c_0^{1, 0}
\\
c_1^{1, 0}
\\
c_2^{1, 0}
\\
d_0^{1, 0}
\\
d_1^{1, 0}
\\
d_2^{1, 0}
\\
\delta
\end{bmatrix}
= B_1 \cdot \vec{v},
\qquad
\begin{bmatrix}
s^{0, 1}
\\
c_0^{0, 1}
\\
c_1^{0, 1}
\\
c_2^{0, 1}
\\
d_0^{0, 1}
\\
d_1^{0, 1}
\\
d_2^{0, 1}
\\
\delta
\end{bmatrix}
= B_2 \cdot \vec{v},
\end{gather*}
where the matrix $B_1$ is given by
\begin{gather*}
B_1 =
\begin{bmatrix}
2 & -1 & 0 & 0 & 0 & -1 & -1 & 1
\\
1 & 0 & 0 & 0 & 0 & -1 & -1 & 1
\\
0 & 0 & 1 & 0 & 0 & 0 & 0 & 0
\\
0 & 0 & 0 & 1 & 0 & 0 & 0 & 0
\\
0 & 0 & 0 & 0 & 1 & 0 & 0 & 0
\\
1 & -1 & 0 & 0 & 0 & 0 & -1 & 1
\\
1 & -1 & 0 & 0 & 0 & -1 & 0 & 1
\\
0 & 0 & 0 & 0 & 0 & 0 & 0 & 1
\end{bmatrix}
,
\end{gather*}
and $B_2$ is given by
\begin{gather*}
%\label{B2}
B_2=
\begin{bmatrix}
2 & -1 & -1 & -1 & 0 & 0 & 0 & 1
\\
1 & 0 & -1 & -1 & 0 & 0 & 0 & 1
\\
1 & -1 & 0 & -1 & 0 & 0 & 0 & 1
\\
1 & -1 & -1 & 0 & 0 & 0 & 0 & 1
\\
0 & 0 & 0 & 0 & 1 & 0 & 0 & 0
\\
0 & 0 & 0 & 0 & 0 & 1 & 0 & 0
\\
0 & 0 & 0 & 0 & 0 & 0 & 1 & 0
\\
0 & 0 & 0 & 0 & 0 & 0 & 0 & 1
\end{bmatrix}
.
\end{gather*}
Note that the matrix $B_2$ is the same as~$B$ in Section~\ref{section3.2}.
Recall that $E^{n, 3}$ series is invariant under the permutations within each sets of parameters $\{c_0, c_1, c_2\}$
and $\{d_0, d_1, d_2\}$.
For $i=0, 1$, set~$s_i$ (resp.~$t_i$) to be the permutation of $c_i$ and $c_{i+1}$ (resp.~$d_i$ and $d_{i+1}$).
The matrix realizations~$S_i$ (resp.~$T_i$) for $s_i$ (resp.~$t_i$) is given by its action on the vector~$\vec{v}$.
For example,
\begin{gather*}
S_0=
\begin{bmatrix}
1 & 0 & 0 & 0 & 0 & 0 & 0 & 0
\\
0 & 0 & 1 & 0 & 0 & 0 & 0 & 0
\\
0 & 1 & 0 & 0 & 0 & 0 & 0 & 0
\\
0 & 0 & 0 & 1 & 0 & 0 & 0 & 0
\\
0 & 0 & 0 & 0 & 1 & 0 & 0 & 0
\\
0 & 0 & 0 & 0 & 0 & 1 & 0 & 0
\\
0 & 0 & 0 & 0 & 0 & 0 & 1 & 0
\\
0 & 0 & 0 & 0 & 0 & 0 & 0 & 1
\end{bmatrix}
\qquad
\text{and}
\qquad
T_0=
\begin{bmatrix}
1 & 0 & 0 & 0 & 0 & 0 & 0 & 0
\\
0 & 1 & 0 & 0 & 0 & 0 & 0 & 0
\\
0 & 0 & 1 & 0 & 0 & 0 & 0 & 0
\\
0 & 0 & 0 & 1 & 0 & 0 & 0 & 0
\\
0 & 0 & 0 & 0 & 0 & 1 & 0 & 0
\\
0 & 0 & 0 & 0 & 1 & 0 & 0 & 0
\\
0 & 0 & 0 & 0 & 0 & 0 & 1 & 0
\\
0 & 0 & 0 & 0 & 0 & 0 & 0 & 1
\end{bmatrix}
.
\end{gather*}
For the action on the variable $x_i$, we set to be $b_2 \cdot x_i = \tilde{x}_i =-m_i \delta-x_i+|M| \delta$ and
identical otherwise.
We denote $I_n$ as the unit $n \times n$ matrix.

We introduce the normalized elliptic hypergeometric series $\widetilde{E}^{n, 3} ((\vec{v}, x))$ as follows
\begin{gather}
\widetilde{E}^{n, 3} \left((\vec{v}, x) \right) := \frac{\prod\limits_{0 \le k \le 2} [[ \delta+s
- c_k]]_{|M|} \left(\prod\limits_{1 \le i \le n} [[ \delta+s-d_k+x_i]]_{m_i} \right)}{\prod\limits_{1 \le i \le n} [[ \delta+s+x_i]]_{m_i}}
\nonumber
\\
\phantom{\widetilde{E}^{n, 3} \left((\vec{v}, x) \right) :=}
\times \ME{n,3}{\{- m_i \delta\}_n\\ \{x_i\}_n}{s}{c_0, c_1, c_2}{d_0,d_1,d_2}.
\label{NormEn}
\end{gather}

Now, we show an invariance property for $E^{n, 3}$ series of type~\eqref{MTEn}.

\begin{proposition}[Hardy type invariant form for $E^{n,3}$ series of type~\eqref{MTEn}]\label{proposition5} Under the balancing condition~\eqref{bc},
$\widetilde{E}^{n, 3} ((\vec{v}, x))$ is invariant under the action of $b_1$, $b_2$, $s_0$, $s_1$, $t_0$ and~$t_1$.
\end{proposition}

In the course of the proof of this proposition, we use the following lemma.
\begin{lemma}\label{lemma6}
If~\eqref{bc} holds, we have the following
\begin{gather}
\label{i1}
1)\quad
[[a]]_{|M|} = (-1)^{|M|} [[ 3s+3 \delta -c_0-c_1-c_2-d_0-d_1-d_2 -a]]_{|M|},
\\
2)\quad
\label{i2}
[[a+x_i]]_{m_i} = (-1)^{m_i} [[ 3s+3 \delta -c_0-c_1-c_2-d_0-d_1-d_2 -a+\tilde{x}_i]]_{m_i}.
\end{gather}
\end{lemma}

\begin{proof}
Since $[[x]]$ is a~odd function of $x$,
\begin{gather*}
[[a]]_{|M|}= [[a]][[a+\delta]] \cdots [[a+(|M|-1) \delta]]
= (-1)^{|M|} [[-a +(1 -|M|) \delta]] \cdots [[-a]]
\\
\phantom{[[a]]_{|M|}}{}
=(-1)^{|M|} [[-a +(1 -|M|) \delta]]_{|M|}.
\end{gather*}
By the balancing condition~\eqref{bc}, we have
\begin{gather*}
-a+(1-|M|) \delta = 3 \delta+3 s -c_0-c_1-c_2-d_0-d_1-d_2-a.
\end{gather*}
Thus we have~\eqref{i1}.
Further, one can check~\eqref{i2} in a~similar fashion.
\end{proof}

\begin{proof}[Proof of Proposition~\ref{proposition5}]
It is not hard to see in the case of $s_0$, $s_1$, $t_0$, $t_1$ since $\widetilde{E}^{n,3}$~\eqref{NormEn} is symmetric with
respect to the subscript~$k$.
For the case of~$b_2$,
\begin{gather*}
\widetilde{E}^{n, 3} \left(b_2 \cdot (\vec{v}, x) \right) = \widetilde{E}^{n, 3} \left((B_2 \vec{v}, \tilde{x})\right)
\\
\qquad
= \frac{\prod\limits_{0 \le k \le 2} [[ \delta+s-c_k]]_{|M|} \left(\prod\limits_{1 \le i \le n}
[[\delta+s+d_k-d_0-d_1-d_2+x_i]]_{m_i} \right)}{\prod\limits_{1 \le i \le n} [[ \delta+s
- d_0-d_1-d_2+ x_i]]_{m_i}}
\\
\qquad
\phantom{=}{}
\times  \ME{n,3}{\{- m_i \delta\}_n
\\
\{\tilde{x}_i\}_n}{s^{0,1}}{c_0^{0,1}, c_1^{0, 1}, c_2^{0, 1}}{d_0^{0, 1}, d_1^{0, 1}, d_2^{0, 1}}
\\
\qquad{}
=\frac{\prod\limits_{0 \le k \le 2} [[\delta+s-c_k]]_{|M|} \left(\prod\limits_{1 \le i \le n}
[[\delta+s+d_k-d_0-d_1-d_2+x_i]]_{m_i} \right)}{ \prod\limits_{1 \le i \le n} [[ \delta+s
- d_0-d_1-d_2+ x_i]]_{m_i}}
\\
\qquad
\phantom{=}{}\times \prod\limits_{1 \le i \le n}
\frac{[[\delta+s-d_0-d_1-d_2 +x_i, \delta+s -d_0 +x_i]]_{m_i}}{[[\delta+s-d_0-d_1 +x_i,\delta+s-d_0-d_2 +x_i]]_{m_i}}
\\
\qquad\phantom{=}{}\times \prod\limits_{1 \le i \le n} \frac{[[\delta+s-d_1 +x_i, \delta+s-d_2+x_i]]_{m_i}}{[[\delta+s-d_1-d_2
+x_i, \delta+s +x_i]]_{m_i}}
\\
\qquad\phantom{=}{}\times  \ME{n,3}{\{- m_i \delta\}_n
\\
\{{x}_i\}_n}{s}{c_0, c_1, c_2}{d_0, d_1, d_2}
=\widetilde{E}^{n, 3} \left((\vec{v}, x) \right).
\end{gather*}
Here we used the transformation~\eqref{KN} and Lemma~\ref{lemma6}.
For $b_1$, one can check similarly.
\end{proof}

Here we shall investigate the compositions of the transformations $b_1$, $b_2$, $s_0$, $s_1$, $t_0$ and $t_1$.

\begin{lemma}
The relations $b_1^2 = b_2^2 = s_0 ^2 = s_1^2 = t_0^2 = t_1^2 = \id$
holds $($where $\id$
stands for identical as a~transformation$)$.
Thus the set of the transformations $\{b_1, b_2, s_0, s_1, t_0, t_1\}$ constitutes the generators of a~Coxeter group
by compositions.
\end{lemma}

\begin{proof}
For $s_0$, $s_1$, $t_0$ and $t_1$, it is obvious since these are the permutation for the coordinates.
For~$b_1$ and~$b_2$, we can check by direct computations of matrices~$B_1$ and~$B_2$ that $B_1^2 = B_2^2= I_8$.
\end{proof}

\begin{remark}
Recall that the variables $x_i$ in Bailey~II~\eqref{KN} change to $\tilde{x_i} = m_i-|M|-x_i$.
It is easy to see that by iterating twice,
\begin{gather*}
\tilde{\tilde{x_i}} = m_i-|M|-\tilde{x_i} = x_i.
\end{gather*}
Thus we see that it turn out to be identity as transformation for $E^{n, 3}$ series by iterating Bailey~II~\eqref{KN}
twice.
\end{remark}

Let $G_r$ be the group generated by $b_1$, $b_2$, $s_0$, $s_1$, $t_0$ and $t_1$.
Now we shall give the relations between the generators of the group $G_r$.
By def\/inition of $s_0$, $s_1$, $t_0$ and $t_1$, the following two braid relations hold:
\begin{gather}
\label{BrRel}
(s_0 s_1)^3 = (t_0 t_1)^3 = \id.
\end{gather}
Note also that, for $i, j \in \{0, 1\}$, $s_i$ and $t_j$ mutually commute.
Other relations, among $b_1$, $b_2$ and others, can be summarized as follows:

\begin{lemma}
The relations
\begin{gather*}
%\label{RelB1}
(b_1 s_0)^3 =(b_1 t_0)^3 = \id
\end{gather*}
hold.
Other pairs of generators of $G_r$ commute.
Especially, $b_2$ commutes with any other ge\-ne\-rators.
\end{lemma}

\begin{proof}
One can check by direct computation for the matrix realization given above.
So we shall leave to readers.
\end{proof}

We def\/ine the mapping $\pi$ as
\begin{gather*}
%\label{TaioPi}
b_1 \mapsto
\sigma_3,
\qquad
 b_2 \mapsto \tau,
\qquad
s_i \mapsto \sigma_{2-i},
\qquad
 t_i \mapsto
\sigma_{4+i},
\qquad
 i= 0, 1 .
\end{gather*}
Then, by braid relations~\eqref{BrRel} and two lemmas above, we see that the following relation holds:
\begin{gather*}
\begin{cases}
\sigma_i \not= \id,
\qquad
\sigma_i^2 = \id, &  i = 1, 2, 3, 4, 5 ,
\\
\sigma_{i}
\sigma_{i+1}
\sigma_{i}
=
\sigma_{i+1}
\sigma_{i}
\sigma_{i+1},
&  i = 1, 2, 3, 4 ,
\\
\sigma_{i}
\sigma_{j}
=
\sigma_{j}
\sigma_{i},
&  |i-j| \ge 2 ,
\end{cases}
\end{gather*}
and $\tau^2 = \id$.
In other words, $\{\sigma_j\}
$ and $\{\tau\}$ is a~realization of $\mathfrak{S}_6$ and $\mathfrak{S}_2$.
Thus we have

\begin{proposition}\label{proposition3.3}
The group $G_r$ is isomorphic to the direct product of $\mathfrak{S}_6$ and $\mathfrak{S}_2$.
\end{proposition}

To summarize the results here, we state the following:

\begin{theorem}
Under the balancing condition~\eqref{bc}, $ {\tilde{E}^{n,3} \left((\vec{v}, x) \right)}$ is invariant under the action of the direct product of $\mathfrak{S}_6$ and
$\mathfrak{S}_2$ realized by the mapping $\pi^{-1}$.
\end{theorem}

We are going to classifying non-trivial transformations for $E^{n,3}$ of rectangular type~\eqref{MTEn} by using the
realization~$\sigma_i$ and~$\tau$.

Proposition~\ref{proposition3.3} tells us that the group $G_r$ of the symmetry of $E^{n, 3}$ of type~\eqref{MTEn} is isomorphic to a~direct
product of the $\mathfrak{S}_6$ and $\mathfrak{S}_2$ and is of order $6 ! \times 2! = 1440$.
Recall again that $E^{n,3}$ series of rectangular type is symmetric with respect to the $c_k$, $k = 0, 1, 2$
and $d_k$, $k = 0, 1, 2$.
Then it is not hard to see that the the right action of $\sigma_1$, $\sigma_2$, $\sigma_4$ and $\sigma_5$ corresponds
to the permutations of the subscript in the sets of parameters
$\{c_0, c_1, c_2\}
$ and $\{d_0, d_1, d_2\}$ and the left action corresponds to the permutation of the location of coordinates.
Thus our problem turns out to give an orbit decomposition of the double coset $H_r \backslash G_r /H_r$, where $H_r$ is
a~subgroup generated by $\sigma_1$, $\sigma_2$, $\sigma_4$, $\sigma_5$,
which is isomorphic to a~direct product of two $\mathfrak{S}_3$.
The representatives of orbits in $H_r \backslash G_r / H_r$ are given by the following:
\begin{gather*}
%\label{RepAnB}
0)
\quad
\omega_0 = \id,
\\
1)
\quad
\omega_1=\sigma_3,
\\
2)
\quad
\omega_2= \sigma_3 \sigma_4 \sigma_2 \omega_1 = \sigma_3 \sigma_4 \sigma_2 \sigma_3,
\\
3)
\quad
\omega_3 = \sigma_3 \sigma_4 \sigma_5 \sigma_2 \sigma_1 \omega_2 = \sigma_3 \sigma_4 \sigma_5 \sigma_2 \sigma_1
\sigma_3 \sigma_4 \sigma_2 \sigma_3.
\end{gather*}
Since the element $\tau$ commutes with $\sigma_i$ for all $i=1,2,3,4,5$,
one can f\/ind that $\tau$ is also commutative with all the representatives $\omega_i$
for $i=0,1,2,3$.

Before going to present a~list of non-trivial possible transformations for $E^{n,3}$ series of ty\-pe~\eqref{MTEn}, we
def\/ine a~transformation $\pi^{r,t}$ as~$\pi^{r, t}: = \tau^t \omega_r$ for $r= 0, 1, 2, 3$ and $t=0, 1$.
We call the transformation formula corresponding to $\pi^{r,t}$ as~$T(r,t)$ and express it as follows
\begin{gather*}
%\label{Exp}
\ME{m,3}{\{- m_i \delta\}_n
\\
\{x_i\}_n}{s}{c_0, c_1, c_2}{d_0,d_1,d_2} =
\\
\qquad{}
P^{r, t} (x; s; C; D)  \ME{m,3}{\{- m_i \delta\}
\\
\{{x}_i^t\}}{{s^{r,t}}}{c_0^{r,t}, c_1^{r, t}, c_2^{r, t}}{d_0^{r,t}, d_1^{r, t}, d_2^{r, t}},
\end{gather*}
where $s^{r, t}$, $c_k^{r, t}$, $d_k^{r, t}$ $(k=0,1,2)$ is parameters associated to the transformation $\pi^{r,t}$ and
$P^{r, t} (x; s;$
$C; D)= P^{r, t} (\{x_i\}_n; s; c_0, c_1, c_2;d_0, d_1, d_2)$ is the corresponding product factor.
For the variab\-les~$x_i^t$, we set $x_i^0= x_i$ and $x_i^1 = \widetilde{x}_i = (|M|- m_i) \delta -x_i$.
Note that Bailey~I~\eqref{MN} and Bailey~II~\eqref{KN} correspond to $T(1,0)$ and $T(0, 1)$ respectively.
Note also that $T(0,0)$ is identical.

Here we exhibit a~list of the product factors and transformations for parameters in $T(r,t)$.
In order to simplify each product factor, we frequently use Lemma~\ref{lemma6}.
Note that the expressions of each product factors have ambiguity because of the balancing condition~\eqref{bc}.

{\bf $\boldsymbol{T(2, 0)}$}

$\bullet$~Product factor
\begin{gather*}
%\label{ProdCB}
P^{2, 0} (x; s; C; D)
= (-1)^{|M|} \frac {[[d_2, \delta+s-c_2-d_0, \delta+s-c_2-d_1]]_{|M|}}{[[\delta+s-c_0, \delta+s-c_1,
\delta+s-c_2]]_{|M|}}
\\
\hphantom{P^{2, 0} (x; s; C; D)=}{} \times \prod\limits_{1 \le i \le n} \left[ \frac {[[\delta+s+x_i, 2\delta+2 s -c_0-d_0-d_1-d_2 +x_i]]_{m_i}}
{[[\delta+s-d_0 +x_i, \delta+s-d_1 +x_i]]_{m_i}} \right.
\\
\hphantom{P^{2, 0} (x; s; C; D)=}{}\times \left.
\frac{[[2 \delta+2 s -c_1 -d_0-d_1-d_2 +x_i]]_{m_i}}{[[3 \delta+3 s-c_0 -c_1 -d_0-d_1-2 d_2 +x_i]]_{m_i}}
\right].
\end{gather*}

$\bullet$~Parameters
\begin{gather*}
\begin{bmatrix}
s^{2, 0}
\\
c_0^{2, 0}
\\
c_1^{2, 0}
\\
c_2^{2, 0}
\\
d_0^{2, 0}
\\
d_1^{2, 0}
\\
d_2^{2, 0}
\end{bmatrix}
=
\begin{bmatrix}
3 s+2 \delta-c_0 -c_1-d_0- d_1-2d_2
\\
s+\delta-d_0-d_2
\\
s+\delta-d_1-d_2
\\
c_2
\\
s+\delta-c_0-d_2
\\
s+\delta-c_1-d_2
\\
2s+2 \delta-c_0-c_1-d_0-d_1-d_2
\end{bmatrix}.
\end{gather*}

{\bf $\boldsymbol{T(3, 0)}$}

$\bullet$~Product factor
\begin{gather*}
%\label{ProdLE}
P^{3, 0} (x; s; C; D)
= (-1)^{|M|} \frac{[[d_0, d_1, d_2]]_{|M|}}{[[\delta+s-c_0, \delta+s-c_1, \delta+s-c_2]]_{|M|}}
\\
\hphantom{P^{3, 0} (x; s; C; D)=}{}\times \prod\limits_{1 \le i \le n} \left[ \frac{[[\delta+s+x_i, 2 \delta+2 s -c_0-d_0-d_1-d_2 +x_i]]_{m_i}}
{[[\delta+ s-d_0 +x_i, \delta+s- d_1 +x_i]]_{m_i}} \right.
\\
\hphantom{P^{3, 0} (x; s; C; D)=}{} \times\! \left.
\frac{[[2 \delta\!+2 s -c_1\! -d_0-d_1\!-d_2\! +x_i, 2 \delta\!+2 s -c_2\!-d_0-d_1\!-d_2\! +x_i]]_{m_i}}
{[[\delta+ s-d_2\! +x_i, 4\delta+4 s -c_0 -c_1\!-c_2\!-2 d_0\!-2 d_1\!-2 d_2\! +x_i]]_{m_i}} \right]\!.
\end{gather*}

$\bullet$~Parameters
\begin{gather*}
\begin{bmatrix}
s^{3, 0}
\\
c_0^{3, 0}
\\
c_1^{3, 0}
\\
c_2^{3, 0}
\\
d_0^{3, 0}
\\
d_1^{3, 0}
\\
d_2^{3, 0}
\end{bmatrix}
=
\begin{bmatrix}
4 s+3 \delta-c_0 -c_1 -c_2-2 d_0- 2 d_1-2d_2
\\
s+\delta-d_0-d_1
\\
s+\delta-d_0-d_2
\\
s+\delta-d_1-d_2
\\
2s+2 \delta-c_0-c_1-d_0-d_1-d_2
\\
2s+2 \delta-c_0-c_2-d_0-d_1-d_2
\\
2s+2 \delta-c_1-c_2-d_0-d_1-d_2
\end{bmatrix}.
\end{gather*}

{\samepage {\bf $\boldsymbol{T(1, 1)}$}

$\bullet$~Product factor
\begin{gather*}
%\label{ProdKNMN}
P^{1, 1} (x; s; C; D)
= \frac {[[\delta+s-c_1-d_0, \delta+s-c_2-d_0]]_{|M|}}{[[\delta+s-c_1, \delta+s-c_2]]_{|M|}}
\\
\hphantom{P^{1, 1} (x; s; C; D)=}{}\times \prod\limits_{1 \le i \le n} \frac{[[\delta+s+x_i, c_0 +x_i ]]_{m_i}}{[[\delta+s-d_0 +x_i, \delta+s-d_1
+x_i]]_{m_i}}
\\
\hphantom{P^{1, 1} (x; s; C; D)=}{}\times \prod\limits_{1 \le i \le n} \frac{[[\delta+s- d_0-d_2 +x_i, \delta+s-d_0-d_1+x_i]]_{m_i}}{[[\delta+s
-d_2 +x_i, c_0- d_0 +x_i]]_{m_i}}.
\end{gather*}}

$\bullet$~Parameters
\begin{gather*}
\begin{bmatrix}
s^{1, 1}
\\
c_0^{1, 1}
\\
c_1^{1, 1}
\\
c_2^{1, 1}
\\
d_0^{1, 1}
\\
d_1^{1, 1}
\\
d_2^{1, 1}
\end{bmatrix}
=
\begin{bmatrix}
3 s+2 \delta-2 c_0 -c_1 -c_2-d_1-d_2
\\
2s+2 \delta -c_0-c_1-c_2-d_1-d_2
\\
s+\delta -c_0-c_2
\\
s+\delta -c_0-c_1
\\
d_0
\\
s+\delta -c_0-d_2
\\
s+\delta -c_0-d_1
\end{bmatrix}.
\end{gather*}

{\bf $\boldsymbol{T(2, 1)}$}

$\bullet$~Product factor
\begin{gather*}
%\label{ProdKNCB}
P^{2, 1} (x; s; C; D)
= (-1)^{|M|} \frac {[[\delta+s-c_2-d_0, d_2, \delta+s-c_2-d_1]]_{|M|}}{[[\delta+s-c_0, \delta+s-c_1,
\delta+s-c_2]]_{|M|}}
\\
\hphantom{P^{2, 1} (x; s; C; D)=}{}\times \prod\limits_{1 \le i \le n} \left[ \frac{[[\delta+s+x_i, c_1 +x_i]]_{m_i}}{[[\delta+s-d_0+x_i, \delta+s
- d_1 +x_i ]]_{m_i}} \right.
\\
\hphantom{P^{2, 1} (x; s; C; D)=}{} \times \left.
\frac {[[c_0 +x_i, \delta+s-d_0-d_1+x_i]]_{m_i}}{[[\delta+s-d_2 +x_i,-\delta-s+c_0+c_1+d_2+x_i]]_{m_i}} \right].
\end{gather*}

$\bullet$~Parameters
\begin{gather*}
\begin{bmatrix}
s^{2, 1}
\\
c_0^{2, 1}
\\
c_1^{2, 1}
\\
c_2^{2, 1}
\\
d_0^{2, 1}
\\
d_1^{2, 1}
\\
d_2^{2, 1}
\end{bmatrix}
=
\begin{bmatrix}
4 s+3 \delta-2 c_0-2 c_1-c_2-d_0- d_1-2d_2
\\
2s+2 \delta-c_0-c_1 -c_2- d_0-d_2
\\
2s+2 \delta-c_0-c_1-c_2- d_1-d_2
\\
s+\delta-c_0-c_1
\\
s+\delta-c_0-d_2
\\
s+\delta-c_1-d_2
\\
2s+2 \delta-c_0-c_1-d_0-d_1-d_2
\end{bmatrix}.
\end{gather*}

{\bf $\boldsymbol{T(3, 1)}$}

$\bullet$
Product factor
\begin{gather*}
%\label{ProdKNLE}
P^{3, 1} (x; s; C; D)
= (-1)^{|M|} \frac {[[d_0, d_1, d_2]]_{|M|}}{[[\delta+s-c_0, \delta+s-c_1, \delta+s-c_2]]_{|M|}}
\\
\hphantom{P^{3, 1} (x; s; C; D)=}{}\times \prod\limits_{1 \le i \le n} \left[ \frac {[[\delta+s+x_i]]_{m_i}}{[[ -2 \delta-2 s+c_0+c_1+c_2+d_0
+ d_1+d_2 +x_i]]_{m_i}} \right.
\\
\hphantom{P^{3, 1} (x; s; C; D)=}{} \times \left.
\frac {[[c_0 +x_i, c_1 +x_i, c_2+x_i]]_{m_i}}{[[\delta+s-d_0 +x_i, \delta+s-d_1 +x_i, \delta+s-d_2
+x_i]]_{m_i}} \right].
\end{gather*}

$\bullet$
Parameters
\begin{gather*}
\begin{bmatrix}
s^{3, 1}
\\
c_0^{3, 1}
\\
c_1^{3, 1}
\\
c_2^{3, 1}
\\
d_0^{3, 1}
\\
d_1^{3, 1}
\\
d_2^{3, 1}
\end{bmatrix}
=
\begin{bmatrix}
5 s+4 \delta-2c_0 -2c_1-2c_2-2 d_0- 2 d_1-2d_2
\\
2s+2 \delta-c_0-c_1-c_2- d_0-d_1
\\
2s+2 \delta-c_0-c_1-c_2- d_0-d_2
\\
2s+2 \delta-c_0-c_1-c_2- d_1-d_2
\\
2s+2 \delta-c_0-c_1-d_0-d_1-d_2
\\
2s+2 \delta-c_0-c_2-d_0-d_1-d_2
\\
2s+2 \delta-c_1-c_2-d_0-d_1-d_2
\end{bmatrix}.
\end{gather*}

Note that by reversing the order of the summation for $E^{n,3}$ series, namely by replacing $\gamma_i \mapsto m_i-\gamma_i$
and simplifying the factors, we also obtain $T(3, 1)$.
Note also that~\eqref{KN} can be obtained by combining $T(3,0)$ and $T(3, 1)$.

\subsection[The case of triangular $E^{n, 3}$ series]{The case of triangular $\boldsymbol{E^{n, 3}}$ series}
\label{section3.4}

Here, we shall discuss triangular case.
That is the case of $E^{n, 3}$ series of the form
\begin{gather}
\label{LTEn}
\ME{n,3}{\{a_i\}_n
\\
\{x_i\}_n}{{s}}{c_0,c_1,c_2}{-N \delta,d_1,d_2},
\end{gather}
which terminates with respect to the length of multi-indices and is provided the balancing condition
\begin{gather}
\label{lbc}
\sum\limits_{1 \le i \le n}
a_i+ c_0+c_1+c_2+d_1+d_2 =(2+N)\delta+3s.
\end{gather}

In this case, we have also obtained the following $A_n$ elliptic Bailey transformation formulas for $E^{n, 3}$ series
for triangular type~\eqref{LTEn} in~\cite{KajiNou}.

{\bf $\boldsymbol{A_n}$ Bailey transformations for $\boldsymbol{E^{n, 3}}$ series of triangular type} (Theorem~4.1
in~\cite{KajiNou}).
Under the balancing condition~\eqref{lbc}, we have two types of $A_n$ Bailey transformation formulas.

 {\bf $\boldsymbol{A_n}$ Bailey~I}
\begin{gather}
\ME{n,3}{\{a_i\}_n\\ \{x_i\}_n}{{s}}{c_0,c_1,c_2}{-N \delta,d_1,d_2}
\nonumber
\\
\qquad
= \frac {[[\delta+\widetilde{s}-c_1, \delta+\widetilde{s}-c_2]]_N}{[[\delta+s-c_1, \delta+s-c_2]]_N} \prod\limits_{1
\le i \le n} \frac {[[\delta+s+x_i, \delta+\widetilde{s}+x_i-a_i]]_N}{[[\delta+s+x_i-a_i, \delta+\widetilde{s}+x_i]]_N}
\nonumber
\\
\qquad
\phantom{=}{}
\times \ME{n,3}{\{a_i\}_n
\\
\{x_i\}_n}{\widetilde{s}}{\widetilde{c}_0,c_1,{c}_2}{- N \delta,\widetilde{d}_1,\widetilde{d}_2},
\label{LMN}
\end{gather}
where
\begin{gather*}
\widetilde{s}=\delta+2 s-c_2-d_0-d_1,
\qquad
\widetilde{c}_0=\delta+s-d_1-d_2,
\\
\widetilde{d}_1=\delta+s-c_0-d_2,
\qquad
\widetilde{d}_2=\delta+s-c_0-d_1.
\end{gather*}

  {\bf $\boldsymbol{A_n}$ Bailey~II}
\begin{gather}
\ME{n,3}{\{a_i\}_n\\ \{x_i\}_n}{s}{c_0,c_1,c_2}{-N \delta,d_1,d_2}
= \prod\limits_{1 \le i \le n} \left[ \frac {[[\delta+s+x_i, \delta+s+x_i-d_1-d_2]]_N}{[[\delta+s+x_i-d_1,
\delta+s+x_i-d_2]]_N} \right.
\nonumber
\\
\qquad
\phantom{=}{}
\times \left.
\frac {[[\delta+s+x_i-a_i-d_1, \delta+s+x_i-a_i-d_2]]_N}{[[\delta+s+x_i-a_i, \delta+s+x_i-a_i-d_1-d_2]]_N} \right]
\nonumber
\\
\qquad
\phantom{=}{}
\times \ME{n,3}{\{a_i\}_n
\\
\{z_i\}_n}{\widetilde{s}}{\widetilde{c}_0,\widetilde{c}_1,\widetilde{c}_2}{-N \delta,d_1,d_2},
\label{LKN}
\end{gather}
where
\begin{gather*}
\widetilde{s}=\delta+2 s-c_0-c_1-c_2,
\qquad
\widetilde{c}_0=\delta+s-c_1-c_2,
\qquad
\widetilde{c}_1=\delta+s-c_0-c_2,
\\
\widetilde{c}_2=\delta+s-c_0-c_1,
\qquad
z_i= a_i-x_i-|a|,
\qquad
 i=1,\ldots,m .
\end{gather*}

Note that, in the case when $n=1, x_1=0$,~\eqref{LMN} and~\eqref{LKN} reduce to the elliptic Bailey transformation
formula~\eqref{EBaileyT1}.

Recall that, though the right hand side in Bailey~I in rectangular case~\eqref{MN} contains two types of $d_j$'s:
$\widetilde{d}_0 =d_0$ f\/ixed and $\widetilde{d}_j = \delta+s -c_0-d_1-d_2+d_j$, $j=1, 2$, it consists of only
$\widetilde{d}_j$ $(j= 1, 2)$ in triangular case~\eqref{LMN}.
Thus we f\/ind that, on the contrast to rectangular case, the element $t_0=\sigma_4$ lacks in this case.
Thus we have:

\begin{proposition}
The group describing the symmetry for the transformations~\eqref{LMN} and~\eqref{LKN} is isomorphic to $\mathfrak{S}_4
\times (\mathfrak{S}_2)^2$.
\end{proposition}

It is not hard to see that the composition of~\eqref{LMN} and~\eqref{LKN} is the only further non-trivial transformation
which can be obtained
\begin{gather}
\ME{n,3}{\{a_i\}_n\\ \{x_i\}_n}{s}{c_0,c_1,c_2}{-N \delta,d_1,d_2}
\nonumber
\\
\qquad
= \frac {[[2 \delta+ {s}-c_0-c_1-d_1-d_2, 2 \delta+ {s}-c_0-c_2-d_1-d_2]]_N}{[[\delta+s-c_1,
\delta+s-c_2]]_N}
\nonumber
\\
\qquad
\phantom{=}
\times \prod\limits_{1 \le i \le n} \left[ \frac{[[\delta+s+x_i, x_i+ c_0]]_N}{[[\delta+s+x_i-a_i, x_i+c_0-a_i]]_N}
\right.
\nonumber
\\
\qquad
\phantom{=}
\times \left.
\frac {[[\delta+s+x_i -d_1-a_i, \delta+ s+x_i -d_2- a_i]]_N}{[[\delta+s+x_i -d_1, \delta+ s+x_i -d_2]]_N}
\right]
\nonumber
\\
\qquad
\phantom{=}
\times \ME{n,3}{\{a_i\}_n\\
\{{z}_i\}_n}{\widehat{s}}{\widehat{c}_0,\widehat{c}_1,\widehat{c}_2}{- N \delta, \widehat{d}_1,\widehat{d}_2},
\label{LMNKN}
\end{gather}
where
\begin{gather*}
\widehat{s}= 2 \delta+3 s-2 c_0-c_1-c_2-d_1-d_2,
\qquad
\widehat{c}_0=2 \delta+2 s- c_0-c_1-c_2-d_1-d_2,
\\
\widehat{c}_1=\delta+s-c_0-c_2,
\qquad
 \widehat{c}_2=\delta+s-c_0-c_1,
\qquad
\widehat{d}_1=\delta+s-c_0-d_2,
\\
 \widehat{d}_2=\delta+s-c_0-d_1,
\qquad
{z}_i=a_i-x_i-|a|,
\qquad
 i=1,\ldots,n .
\end{gather*}

To simplify the product factor, we used the following lemma which can be proved just in the same line as in the
rectangular case.

\begin{lemma}
If the balancing condition~\eqref{lbc} holds, then we have
\begin{gather*}
[[b]]_N = (-1)^N
[[3 \delta + 3 s - b - (c_0 +c_1 +c_2) -
(d_1 + d_2) - |a|]] _N.
\end{gather*}
\end{lemma}

\subsection{Remarks on results of Section~\ref{section3}}

We close this paper to give some remarks.

\begin{remark}
The transformation $T(1,1)$ in Section~\ref{section3.3} has appeared
as Corollary~4.3 in Rosengren~\cite{RoseKaji} with a~dif\/ferent
expression and the transformation~\eqref{LMNKN} has appeared as Corollary~4.2 in~\cite{RoseKaji}.
\end{remark}

\begin{remark}[in the case when ${n=1}$, ${x_1=0}$]
In this case, $T({2,0})$ and $T({1,1})$ in Section~\ref{section3.3}.
and~\eqref{LMNKN} in Section~\ref{section3.4}.
reduce to~\eqref{mn1EBDT1} in Section~\ref{section3.2}.
$T({3, 0})$ and $T({2, 1})$ reduce to~\eqref{BaileyT3}.
Finally, $T(3,1)$ reduces to~\eqref{BaileyT4}.
Notice that~\eqref{BaileyT4} can also be obtained by reversing order of the summation in the ${}_{10} E_9$ series.
\end{remark}

\begin{remark}[correspondence of the group ${G_r}$ in Section~\ref{section3.3}
and ${G_1}$ in Section~\ref{section3.2}]
By direct computation using the matrix
realization in this paper, one f\/inds that $b_1= \pi^{-1}(\sigma_3)$ can be expressed as
\begin{gather}
\label{hrD5}
\nu^{-1} w_2 \nu^{},
\qquad
\nu = w_4 w_5 w_6 w_3 w_4 w_5.
\end{gather}
The correspondence between the generators of the group $G_r$ in Section~\ref{section3.3}.
and the elements of the group $G_1$ is summarized as follows:
\begin{gather*}
\begin{array}{@{}ccc}
  G_1\simeq W(E_6) &&
G_r \simeq \mathfrak{S}_6 \times \mathfrak{S}_2
\\
  w_1
&
 \longleftrightarrow &
\sigma_2,
\\
  w_3
&
 \longleftrightarrow
&
\sigma_1,
\\
  w_5
&
 \longleftrightarrow
&
\sigma_4,
\\
  w_6
&
 \longleftrightarrow
&
\sigma_5,
\\
  w_2
&
 \longleftrightarrow
&
\tau,
\\
  \nu^{-1}
w_2 \nu
&
 \longleftrightarrow
&
\sigma_3.
\end{array}
\end{gather*}
The correspondence is described diagrammatically as follows:

\begin{center}
\setlength{\unitlength}{0.1mm}
\begin{picture}
(1200, 1000) \put(100, 950){$n = 1$} \put(600, 750){\circle{14}} \put(200, 750){\circle{14}} \put(400, 750){\circle{14}}
\put(800, 750){\circle{14}} \put(1000, 750){\circle{14}} \put(600, 950){\circle{14}} \put(207, 750){\line(1,0){186}}
\put(407, 750){\line(1,0){186}} \put(607, 750){\line(1,0){186}} \put(807, 750){\line(1,0){186}} \put(600,
757){\line(0,1){186}} \put(190, 720){$w_1$} \put(390, 720){$w_3$} \put(590, 720){$w_4$} \put(790, 720){$w_5$} \put(990,
720){$w_6$} \put(610, 940){$w_2$}  \put(580, 600){\huge $\Updownarrow $}   \put(100, 450){$ n \ge 2$} \put(200,
250){\circle{14}} \put(400, 250){\circle{14}} \put(800, 250){\circle{14}} \put(1000, 250){\circle{14}} \put(600,
450){\circle{14}} \put(207, 250){\line(1,0){186}}   \put(807, 250){\line(1,0){186}}  \put(140, 280){$w_1 \leftrightarrow
\sigma_2$}
\put(340, 280){$w_3 \leftrightarrow
\sigma_1$}
\put(740, 280){$w_5 \leftrightarrow
\sigma_4$}
\put(940, 280){$w_6 \leftrightarrow
\sigma_5$}
\put(620, 440){$w_2 \leftrightarrow \tau$} \put(600, 50){\circle{14}} \put(605, 52){\line(1,1){190}} \put(592,
52){\line(-2,1){385}} \put(620, 30){$\nu^{-1} w_2 \nu \leftrightarrow
\sigma_3$}
\end{picture}

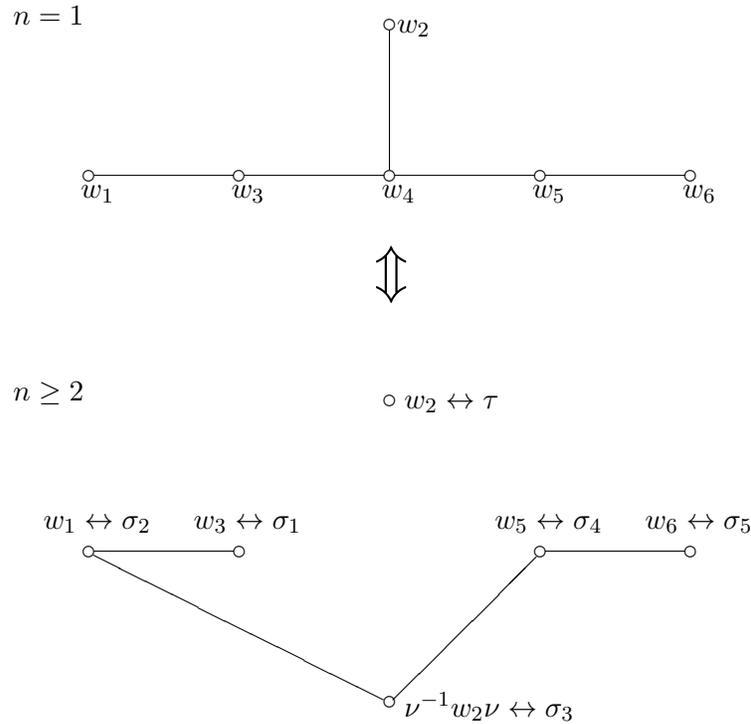
{\small Correspondence of the elements}
\end{center}

Thus we f\/ind the group generated by $w_1$, $w_3$, $w_4$, $w_5$ and $\nu^{-1}\sigma_2 \nu^{}$
is isomorphic to the symmetric group $\mathfrak{S}_6$ and all the generators commute with~$w_2$.

Furthermore, $r_\beta = \nu^{-1} w_2 \nu \in W(E_6)$~\eqref{hrD5} is the ref\/lection of the root $\beta= \alpha_2+\alpha_3+2 \alpha_4+2 \alpha_5+\alpha_6$.
Note that $\beta$ is the highest root of the root system $D_5$ whose roots are $\alpha_2$, $\alpha_3$, $\alpha_4$,
$\alpha_5$ and~$\alpha_6$:
\begin{center}
\setlength{\unitlength}
{0.1mm}
\begin{picture}
(1100, 400) \put(600, 200){\circle{14}} \put(200, 333){\circle{14}} \put(400, 200){\circle{14}} \put(800,
333){\circle{14}} \put(800, 67){\circle{14}} \put(200, 67){\circle{14}} \put(207, 330){\line(3,-2){186}} \put(407,
200){\line(1,0){186}} \put(607, 203){\line(3, 2){186}} \put(607, 197){\line(3, -2){186}} \put(207, 70){\line(3, 2){186}}
\put(190, 303){$w_6$} \put(390, 170){$w_5$} \put(590, 170){$w_4$} \put(790, 37){$w_3$} \put(790, 303){$w_2$} \put(190,
37){$r_\beta$}
\end{picture}

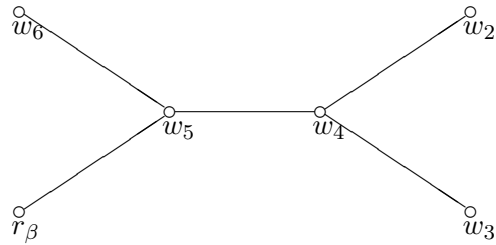
{\small Extended Dynkin diagram of $D_5$}
\end{center}
Note also that the expression in~\eqref{hrD5} of $r_\beta$ is reduced.
\end{remark}

\subsection*{Acknowledgments}

I would like to express my sincere thanks to Professors David Bessis, Christian Krattenthaler, Masato Okado and Hiroyuki
Yamane and, in particular, Professor Kenji Iohara for crucial comments and fruitful discussions on some Coxeter groups
that appear in this paper.
I also thank to anonymous referees to pointing out the errors of the previous version of this paper and useful
suggestions to improve the descriptions.

\pdfbookmark[1]{References}{ref}
\LastPageEnding

\end{document}